\documentclass[a4]{amsart}
\usepackage{amssymb}
\usepackage{amscd}
\usepackage{verbatim,ifthen}
\usepackage{graphicx}
\usepackage{xypic}
\usepackage{cite}
\usepackage{amsthm}

\usepackage{latexsym}

\newtheorem*{lemma*}{Lemma}

\DeclareMathOperator*{\esssup}{ess\,sup}

\newtheorem{thm}{Theorem}

\newtheorem{lemma}[thm]{Lemma}
\newtheorem{cor}[thm]{Corollary}
\newtheorem{rem}[thm]{Remark}
\theoremstyle{definition}
\newtheorem{defn}[thm]{Definition}

\usepackage{hyperref}

\begin{document}

\title[A Weyl Pseudodifferential Calculus associated with exponential weights on $\mathbb{R}^d$]{A Weyl Pseudodifferential Calculus associated with exponential weights on $\mathbb{R}^d$} 

\author{Sean Harris}
\address{Hanna Neumann Building \#145, Science Road
The Australian National University
Canberra ACT 2601.}
\email{Sean.Harris@anu.edu.au}     

\date{\today}

\thanks{The author gratefully acknowledges financial support by the discovery Grant DP160100941 of the Australian Research Council. This research is also supported by an Australian Government Research Training Program (RTP) Scholarship.}

 \begin{abstract}
We construct a Weyl pseudodifferential calculus tailored to studying boundedness of operators on weighted $L^p$ spaces over $\mathbb{R}^d$ with weights of the form $\exp(-\phi(x))$, for $\phi$ a $C^2$ function, a setting in which the operator associated to the weighted Dirichlet form typically has only holomorphic functional calculus. A symbol class giving rise to bounded operators on $L^p$ is determined, and its properties analysed. This theory is used to calculate an upper bounded on the $H^\infty$ angle of relevant operators, and deduces known optimal results in some cases. Finally, the symbol class is enriched and studied under an algebraic viewpoint.
\end{abstract}

 
 \maketitle

\tableofcontents

\section{Introduction}

We construct a Weyl pseudodifferential calculus on typically non-doubling measure spaces over $\mathbb{R}^d$, with an aim to study the $L^p$ behaviour of the natural analogue of the Laplacian in such contexts. Typically, this analogue of the Laplacian is such that bounded spectral multipliers on $L^p$  have to be of holomorphic type for $p \neq 2$. From a PDE point of view, the operators we consider are perturbations of the Laplacian by an unbounded drift term, including the classical finite-dimensional Ornstein-Uhlenbeck operator.

In the paper \cite{vNPWeylGaussian}, van Neerven and Portal introduce a Gaussian Weyl calculus to study the classical Ornstein-Uhlenbeck operator $L=-\Delta + x \cdot \nabla$ on $L^p\left(\mathbb{R}^d, (2\pi)^{-\frac{d}{2}}\exp\left(-\frac{x^2}{2}\right)dx\right)$. This approach retrieves important known results about the Ornstein-Uhlenbeck semigroup, such as boundedness $L^p \to L^q$ and optimal domains of holomorphy for the semigroup, using analytic arguments as simple as Schur estimates. With classical approaches, such properties are difficult to prove.

The classical Ornstein-Uhlenbeck operator can be written in terms of a pair of operators P (momentum) and Q (position), satisfying the Heisenberg commutation relations. The Weyl calculus examined by van Neerven and Portal in \cite{vNPWeylGaussian} is a certain choice of joint functional calculus for this pair $(Q,P)$, that is, a way to assign to a suitable function $a:\mathbb{R}^{2d} \to \mathbb{C}$ a bounded operator $a(Q,P)$. It was their philosophy that studying $L$ via studying the joint functional calculus was more natural, as it separates the strong algebraic properties of the pair $(Q,P)$ from the analytic issues found in proving properties of $L$ directly. Essentially, studying $L$ directly is forgetting that it has its roots in an algebraically rich setting. The approach used by van Neerven and Portal in \cite{vNPWeylGaussian} is well-adapted to studying the standard Ornstein-Uhlenbeck semigroup, thanks to an exact formula for the semigroup in the $(Q,P)$-calculus.

To generalise the theory developed in \cite{vNPWeylGaussian}, van Neerven and Portal consider in \cite{vNPWeylPairs} pairs of $d$-tuples of operators on a Banach space, which generate uniformly bounded groups satisfying certain commutation relations, which they refer to as a Weyl pair. This path of generalisation has been fruitful, and proves that the sum of squares of all $2d$ operators comprising a Weyl pair with a bounded Weyl calculus has a bounded H\"{o}rmander functional calculus (under some mild conditions).

However, in the Gaussian case, the natural position and momentum operators generate uniformly bounded groups on $L^p$ if and only if $p=2$. In fact, $\exp(i\xi P)$ is bounded on $L^p$ for $p\neq 2$ if and only if $\xi = 0$, in which case $\exp(i 0 P)$ is the identity. Thus the theory developed in \cite{vNPWeylPairs} cannot possibly be applied in Gaussian situations. In this paper we consider such cases.

We will work on measure spaces of the form $(\mathbb{R}^d, \mathcal{B}, \mu)$, where $\mathcal{B}$ is the Borel $\sigma$-algebra and $d\mu = \exp(-\phi(x))dx$ for $\phi \in C^2(\mathbb{R}^d)$ approximately quadratic (see Remark \ref{remAffine}).  We introduce a generalised Weyl pair $(Q,P)$ associated with such a measure via a specific unitary equivalence on $L^2(\mu)$ to the standard Weyl pair on $\mathbb{R}^d$ equipped with Lebesgue measure. Typically, the pair $(Q,P)$ will not generate uniformly bounded groups on $L^p(\mu)$ for $p \neq 2$. Thus, our generalised Weyl calculus will be developed as an extension theory: we work on $L^p(\mu)\cap L^2(\mu)$ and find conditions under which we have a bounded extension to $L^p(\mu)$.

For a function $M:\mathbb{R}^d \to \mathbb{R}^d$, we introduce the normed vector space $HS_0(M) \subset L^\infty\left(\mathbb{R}^{2d}\right)$ of Holomorphic Strip symbols. The main theorem of this paper (Theorem \ref{thmBoundedness0}) proves that for a correctly chosen function $M$ based on $p$ and $\phi$, $a \in HS_0(M)$ implies that $a(Q,P)$ extends from $L^p(\mu)\cap L^2(\mu)$ to a bounded operator on $L^p(\mu)$, with norm bounded by $||a||_{HS_0(M)}$. We prove that symbols in $HS_0(M)$ have a certain holomorphic extendability property, which is reminiscent of the optimal functional calculus result for the classical Ornstein-Uhlenbeck operator in the Gaussian case. We also deduce a simple condition on a set $A \subset HS_0(M)$ which implies that $\{a(Q,P); a \in A\}$ is R-bounded on $L^p(\mu)$.

This theory is used to give a short proof of an upper bound on the optimal angle of the bounded functional calculus result for the Ornstein-Uhlenbeck operator associated with $\phi$ of the form $\phi(x) = xN(x) + lx$, where $N:\mathbb{R}^d \to \mathbb{R}^d$ is a positive definite real-symmetric linear operator, and $l \in \mathbb{R}^d$. This result returns the known optimal angle for the classical Ornstein-Uhlenbeck operator, for which $N=\frac{1}{2}I, l=0$. The proof depends on knowing the Weyl symbol for the semigroup generated by the Ornstein-Uhlenbeck operator, which we provide. This explicit knowledge of the symbol in the classical case is also crucial to the results of the paper \cite{vNPWeylGaussian}.

The space $HS_0(M)$ is then examined further, with an aim to find symbols for the semigroup generated by the relevant Ornstein-Uhlenbeck operator outside of the quadratic case. We show that $HS_0(M)$ can be expanded and endowed with a certain product to form the unital Banach algebra $(HS(M), \#)$, with $HS(M) \subset L^\infty\left(\mathbb{R}^{2d}\right)$ and $\#$ the Moyal product, in which case the generalised Weyl calculus can be extended in a natural way such that it becomes a contractive Banach algebra homomorphism $HS(M) \to B(L^p(\mu))$.

Some ideas are then presented, about how the symbol class $HS(M)$ and its properties as a Banach algebra of functions could be used to study Ornstein-Uhlenbeck like operators, such as determining symbols for relevant semigroups and a factorisation of the functional calculus through the space $HS(M)$.

\section{Initial Definitions}\label{secInitialDefinitions}

The work of van Neerven and Portal in \cite{vNPWeylPairs} develops a Weyl calculus for $2d$-tuples of operators satisfying the following definition:

\begin{defn}\label{defnWeylPair}
A $2d$-tuple $(Q,P) = (Q_1, \ldots, Q_d, P_1, \ldots, P_d)$ of closed, densely defined operators on a Banach space is called a Weyl pair if each operator generates a uniformly bounded group $\exp(ix_iQ_i), \exp(i\xi_iP_i)$ which satisfy the integrated canonical commutation relations
\[\begin{array}{c}
\exp(ix_iQ_i)\exp(ix_jQ_j) = \exp(ix_jQ_j)\exp(ix_iQ_i), \forall i,j = 1,... d\\
\exp(i\xi_iP_i)\exp(i\xi_jP_j) = \exp(i\xi_jP_j)\exp(i\xi_iP_i), \forall i,j = 1,... d\\
\exp(ixQ)\exp(i\xi P) = \exp(ix\xi)\exp(i\xi P)\exp(ixQ).\\
\end{array}\]
In which case we define
\[\exp(i(xQ + \xi P)) = \exp(\frac{1}{2}ix\xi)\exp(ixQ)\exp(i\xi P), \forall x,\xi \in \mathbb{R}^d.\]
\end{defn}

This definition of $\exp(i(xQ + \xi P))$ for Weyl pairs can be motivated by the Baker-Campbell-Hausdorff formula, for example, noting that formally differentiating the integrated commutation relations produces the formal commutation relation ``$[Q_i, P_j] = i \delta_{ij}I$". It is easy to check that this definition does make the set of operators $\{\exp(i(xQ + \xi P)); x, \xi \in \mathbb{R}^d\}\cup\{ \lambda I; \lambda \in \mathbb{C}, |\lambda|=1\}$ into a (non-commutative) group, in fact a representation of the $(2d+1)$-dimensional Heisenberg group. From here, the Weyl calculus is defined as:

\begin{defn}\label{defnWeylCalc}
For a Weyl pair $(Q,P)$ on a Banach space, we define the bounded operator $a(Q,P)$ for $a \in \mathcal{S}(\mathbb{R}^{2d})$ via the formula
\[a(Q,P) = \frac{1}{(2\pi)^{d}}\int_{\mathbb{R}^{2d}} \mathcal{F}a(x, \xi) \exp(i(xQ + \xi P)) dxd\xi,\]
where $\mathcal{F}$ denotes the Fourier transform, normalised such that equality holds above with  $Q$ and $P$ replaced by elements of $\mathbb{R}^d$.
The map $a \mapsto a(Q,P)$ is called the Weyl pseudodifferential calculus, or the joint functional calculus, for the Weyl pair $(Q,P)$.
\end{defn}

That this definition makes sense follows from boundedness of the group $\exp(i(xQ + \xi P))$ and integrability of the Schwartz function $\mathcal{F}a$. As an example of a Weyl pair, consider the position and momentum operators $(X,D)$ on the Euclidean Lebesgue space $L^p(\lambda)$ ($p \in [1, \infty]$) given by $X_if(x) = x_if(x), P_if(x) = -i\frac{\partial}{\partial x_i}f(x)$ equipped with their natural domains. These generate the groups of phase shift and translation respectively, both of which are bounded. That they satisfy the integrated canonical commutation relations is a simple exercise, which once checked shows that $(X,D)$ on $L^p(\lambda)$ are a Weyl pair, for any value of $p \in [1, \infty]$. The pair $(X,D)$ will be referred to as the standard pair, and their Weyl calculus is the standard Weyl calculus (see for example \cite{HormanderWeylCalc}).

The body of the work of [vN,P] follows from the next theorem, displaying the semigroup generated by $\frac{1}{2}(Q^2+P^2-d)$ in the joint functional calculus, which turns out to be the correct expression for the Ornstein-Uhlenbeck operator in the case they consider. This formula has been known by physicists for the standard pair $(Q,P)=(X,D)$ (in which case $\frac{1}{2}(X^2+D^2-d)$ is the harmonic oscillator operator minus $d/2$ times the identity), and relies heavily on the simple algebraic nature of $\frac{1}{2}(Q^2+P^2)$ in terms of $(Q,P)$.

\begin{thm}\label{TheoremSGEasy}
For $t > 0$, let $a_t:\mathbb{R}^{2d} \to \mathbb{C}$ be the function $(x, \xi) \mapsto \left(\frac{2}{1+e^{-t}}\right)^d \exp(- s_t (x^2 + \xi^2))$, where $s_t = \frac{1-e^{-t}}{1+e^{-t}}$. Then $\{a_t(Q,P)\}_{t > 0}$ is the semigroup generated by $\frac{1}{2}(Q^2+P^2 - d)$, on its natural domain inherited from the domains of $Q$ and $P$.
\end{thm}

In this paper we wish to develop a different side to the Weyl calculus, in which we do not require the bounded group assumption of Definition \ref{defnWeylPair}. The main reason for this is because the key example we wish to apply the Weyl calculus to - the original motivation of study for van Neerven and Portal in \cite{vNPWeylGaussian} - is the natural position and momentum operator pair in the standard Gaussian weighted spaces $L^p\left(\mathbb{R}^d, \exp\left(-\frac{x^2}{2}\right)dx\right)$, in which the momentum operator $P$ does not satisfy the bounded group generation property unless $p=2$. However, in this case we can define the Weyl calculus on $L^2$ via Definition \ref{defnWeylCalc} and then determine conditions under which $a(Q,P)$ extends to a bounded operator on $L^p$. The cases we will consider will be based over $\mathbb{R}^d$ equipped with Borel measures of the following form.

\begin{defn}
A function $\phi \in C^2(\mathbb{R}^d)$ will be referred to as a potential. associated with a potential $\phi$ is the Borel measure $\mu$ on $\mathbb{R}^d$ with $d\mu = \exp(-\phi(x))dx$.
\end{defn}

We will generally assume $\phi$ to be a twice differentiable function throughout the paper. Although the initial definitions work for any $C^2$ function $\phi$, to obtain boundedness of operators we will later have to restrict to $\phi$ which is approximately quadratic (see Remark \ref{remAffine}). To be able to relate things to the standard pair, we proceed via unitary equivalence:

\begin{defn}\label{defnU2}
For $c>0$, define the (multiples of) unitary transformations
\[\tilde{U}_2 = A \circ E: L^2(\mu) \to L^2(\lambda)\]
where
\[ \begin{array}{c}
E: L^2(\mu) \to L^2(\lambda)\\
f \mapsto \left(x \mapsto f(x) \exp(-\frac{\phi(x)}{2})\right)\\
\\
A: L^2(\lambda) \to L^2(\lambda)\\
f \mapsto \left(x \mapsto f(cx)\right)\\
\end{array}\]
\end{defn}

\begin{defn}\label{defnUp}
Define the isometry $U_p: L^p(\mu) \to L^p(\lambda)$ by
\[f \mapsto \left(x \mapsto f(x) \exp(-\frac{\phi(x)}{p})\right)\]
\end{defn}

The generalised Weyl pairs we will consider in this context will be
\[\begin{array}{c}
Q_i = \tilde{U}_2^{-1} \circ X_i \circ \tilde{U}_2\\
P_i = \tilde{U}_2^{-1} \circ D_i \circ \tilde{U}_2
\end{array}\]
with domains $\tilde{U}_2^{-1}D(X_i)$, $\tilde{U}_2^{-1}D(D_i)$ respectively. That $(Q,P)$ satisfy the requirements of Definition \ref{defnWeylPair} on $L^2(\mu)$ follows from their unitary equivalence to the standard pair. It is a simple exercise to determine the action of $Q$ and $P$ on their domains. They act as 
\[\begin{array}{c}
(Q_i f)(x) = \frac{x_i}{c}f(x)\\
(P_i f)(x) = -ic\left(\frac{\partial}{\partial x_i} - \frac{1}{2}\frac{\partial \phi}{\partial x_i}\right)f(x).
\end{array}\]
The integral kernel of operators $a(Q,P)$ can be calculated explicitly via unitary equivalence to the standard Weyl Calculus (see, for example, \cite{HormanderWeylCalc}), leading to the formula:
\begin{thm}\label{thmKernel}
For $a \in \mathcal{S}\left(\mathbb{R}^{2d}\right)$, $Q, P$ as above and $f \in L^2(\mu)$, we have for all $y \in \mathbb{R}^d$
\[(a(Q,P)f)(y) = \frac{1}{(2\pi)^{d}c^d} \int_{\mathbb{R}^{2d}}a\left(\frac{x+y}{2c}, \xi\right) \exp\left(-i\xi \left(\frac{x-y}{c}\right)\right) \exp\left(\frac{1}{2}(\phi(y)+\phi(x))\right)f(x) d\xi d\mu(x).\]
\end{thm}

\begin{rem}
The inclusion of $c$ is to satisfy physicists. What many physicists would consider as THE Ornstein-Uhlenbeck operator, and what we shall refer to as the classical Ornstein-Uhlenbeck operator, is  related to the choice $\phi(x) = \frac{x^2}{2}$. In its relationship to Fock spaces, a scaling is introduced to make the classical Ornstein-Uhlenbeck operator appear more symmetric in some sense, which corresponds to taking $c=\sqrt{2}$. However, as we shall see later (take for example Definition \ref{defnMultiplierExt}), in many ways it seems more natural to take $c=1$. We will not include the subscript $c$ as part of the notation, but it will always be lurking in the background ready to be set to $1$ or $\sqrt{2}$. Note that the un-tilde'd $U_2$ falls under Definition \ref{defnUp}, and does not depend on $c$.
\end{rem}

To ensure technical correctness, we need some information about the domains of important operators. The next theorem provides a suitable $p$-independent core for our functional calculus

\begin{thm}\label{thmCore}
The space $\mathcal{C}_\phi = U_2^{-1} C_c^\infty(\mathbb{R}^d)$ is dense in $L^p(\mu)$ for all $p \in [1, \infty)$, and elements of $\mathcal{C}_\phi$ are in $C^2_c(\mathbb{R}^d)$.
\end{thm}

\begin{proof}
Regularity follows from the chain rule and the regularity of $\phi$. To show density in $L^p(\mu)$ we will show density of $U_p\mathcal{C}_\phi$ in $L^p(\lambda)$, employing the isometry of Definition \ref{defnUp}. Note that $U_p\mathcal{C}_\phi$ contains every function of the form $g \exp\left(\left(\frac{1}{2}-\frac{1}{p}\right)\phi \right)$ for $g \in C_c^\infty(\mathbb{R}^d)$, which are all bounded since $\exp\left(\left(\frac{1}{2}-\frac{1}{p}\right)\phi \right)$ is continuous and hence bounded on the compact set $\text{supp}(g)$.
Let $X \subset L^p(\lambda)$ be those functions which have compact support. Note that $X$ is dense in $L^p(\lambda)$ and $U_p\mathcal{C}_\phi$ is contained in $X$, so if we can show $U_p\mathcal{C}_\phi$ is dense in $X$ then we are done. To this effect, fix some $f \in X$. Since $f$ is compactly supported, $\exp\left(-\left|\frac{1}{2}-\frac{1}{p}\right)\phi\right)$ attains a maximum on $\text{supp}(f)$, and so $F(\cdot) = f(\cdot)\exp\left(-\left(\frac{1}{2}-\frac{1}{p}\right)\phi(\cdot)\right) \in L^p(\lambda)$. Since $C_c^\infty(\mathbb{R}^d)$ is dense in $L^p(\lambda)$, we can choose a sequence $\{g_n\} \subset C_c^\infty(\mathbb{R}^d)$ (furthermore, each with support $\text{supp}(g_n) \subset 2\text{supp}(f)$, say) such that $||F-g_n||_{L^p(\lambda)}\to 0$. Thus we have:
\begin{align*}
||f(\cdot)-g_n(\cdot)\exp\left(\left(\frac{1}{2}-\frac{1}{p}\right)\phi(\cdot)\right)||_{L^p(\lambda)}^p &= \int_{2\text{supp}(f)}\left|f(x)-g_n(x)\exp\left(\left(\frac{1}{2}-\frac{1}{p}\right)\phi(x)\right)\right|^p dx \\
&= \int_{2\text{supp}(f)}\exp\left(\left(\frac{p}{2}-1\right)\phi(x)\right)|F(x)-g_n(x)|^p dx \\
&\leq C \int_{2\text{supp}(f)}|F(x)-g_n(x)|^p dx \\
&\leq C ||F - g_n||_{L^p(\lambda)}^p \\
&\to 0.
\end{align*}
So we are done.
\end{proof}

The formula in Theorem \ref{thmKernel} and the formula for the semigroup generated by $\frac{1}{2}(Q^2+P^2-d)$ (Theorem \ref{TheoremSGEasy}) is what allowed van Neerven and Portal to deduce such significant results for the classical Ornstein-Uhlenbeck semigroup via the Weyl calculus in \cite{vNPWeylPairs}. The explicit formula for kernels of $a(Q,P)$ will allow us to work on boundedness on $L^p(\mu)$ of more general symbols.

We will now define the Ornstein-Uhlenbeck operator in our setting. We consider the Dirichlet form $E(f,g) = \int_{\mathbb{R}^d} \nabla f \nabla \overline{g} d\mu$ with domain $\mathcal{C}_\phi$. It follows from the general theory of Dirichlet forms (see \cite{FukushimaOshimaTakeda}) that the operator on $L^2(\mu)$ associated with the form $E$ is positive and generates a positive contraction $C_0$-semigroup, which extends to a positive contraction $C_0$-semigroup on $L^p(\mu)$ for all $p \in (1, \infty)$. We call this semigroup the Ornstein-Uhlenbeck semigroup and denote it by $\exp(-tL)$, and its generator the Ornstein-Uhlenbeck operator which we will denote by $L$. It is a simple computation to find that on $f \in \mathcal{C}_\phi$, $L$ has action
\begin{equation}\label{equationAction}
Lf (x) = -\Delta f(x) + \nabla\phi(x) \cdot\nabla f(x).
\end{equation}
Due to our set-up, we have the fact:
\begin{cor}\label{corSomeHinfty}
For any potential $\phi \in C^2$ and $p \in (1, \infty)$, the corresponding Ornstein-Uhlenbeck operator has a bounded $H^\infty$ functional calculus of some angle less than $\pi/2$.
\end{cor}
This follows from Theorem 10.7.13 of \cite{AnalysisInBanachSpacesV2}, as below.
\begin{thm}\label{thmPosConImpliesHInfty}
Suppose $(\Omega, m)$ is a measure space ($\sigma$-algebra omitted). If an unbounded operator $T$ on $L^p(\Omega, m)$, $p \in (1, \infty)$ generates a positive contraction semigroup, then $T$ has a bounded $H^\infty$ functional calculus of some angle less than $\pi/2$.
\end{thm}

After developing the generalised Weyl calculus associated with the generalised Weyl pair $(Q,P)$ as defined above, we will aim to study the Ornstein-Uhlenbeck operator via the generalised Weyl calculus. In \cite{vNPWeylGaussian}, the formal expression for the classical Ornstein-Uhlenbeck operator $L= \frac{1}{2}(Q^2+P^2-d)$ was very important. This is no longer true in our case. We have as a replacement the following theorem, representing $L$ associated with $\phi$ in the generalised Weyl calculus of the pair $(Q,P)$ associated with $\phi$ (at least formally).

\begin{thm}\label{thmH}
Take $Q,P$ as above, and let $h:\mathbb{R}^{2d} \to \mathbb{R}$ be the function taking value $h(x, \xi) = \frac{\xi^2}{c^2} + C_\phi (x)$, where $C_\phi(x) = \frac{1}{4}|\nabla\phi(cx)|^2 - \frac{1}{2}\Delta \phi(cx)$. For $f \in \mathcal{C}_\phi$, define $h(Q,P)f \in L^1_{loc}(\mathbb{R}^d)$ by $\left((x,\xi) \mapsto C_\phi(x)\right)(Q,P)$ interpreted as a multiplication operator by $C_\phi(x/c)$ and $\left((x,\xi) \mapsto \frac{\xi^2}{c^2}\right)(Q,P)$ interpreted as $\frac{1}{c^2}P^2$. Then $Lf = h(Q,P)f$.
\end{thm}

\begin{proof}
Fix $f \in \mathcal{C}_\phi$. Let $H = U_2 L U_2^{-1}$. Note that $U_2 f =:g \in C_c^\infty(\mathbb{R}^d)$, so we want to consider $H$ acting on $g \in C_c^\infty(\mathbb{R}^d)$. Calculating $(Hg)(x)$ for all $x\in \mathbb{R}^d$ gives:
\[(U_2^{-1} g)(x) = g(x/c)\exp\left(\frac{\phi(x)}{2}\right).\]
So
\begin{align*}
(LU_2^{-1}g)(x) &= \sum_{k=1}^d \left(\left(-\left(\partial_k\right)^2 + \partial_k \phi(\cdot) \partial_k \right)g(\cdot/c)\exp\left(\frac{\phi(\cdot)}{2}\right)\right)(x)\\
&= \sum_{k=1}^d \left(-\partial_k\left( \frac{1}{c}\partial_k g(\cdot/c)\exp\left(\frac{\phi(\cdot)}{2}\right) + \frac{1}{2}g(\cdot/c)\exp\left(\frac{\phi(\cdot)}{2}\right)\partial_k\phi(\cdot) \right)\right)(x) \\
&+ \partial_k \phi(x) \left( \frac{1}{c}\partial_k g(x/c)\exp\left(\frac{\phi(x)}{2}\right) + \frac{1}{2}g(x/c)\exp\left(\frac{\phi(x)}{2}\right)\partial_k\phi(x) \right)\\
&= \sum_{k=1}^d  -\frac{1}{c^2}\partial^2_k g(x/c)\exp\left(\frac{\phi(x)}{2}\right) -  \frac{1}{2c}\partial_k g(x/c)\exp\left(\frac{\phi(x)}{2}\right)\partial_k\phi(x) \\
&- \frac{1}{2c}\partial_k g(x/c)\exp\left(\frac{\phi(x)}{2}\right)\partial_k\phi(x) - \frac{1}{4}g(x/c)\exp\left(\frac{\phi(x)}{2}\right)(\partial_k\phi(x))^2 - \frac{1}{2}g(x/c)\exp\left(\frac{\phi(x)}{2}\right)\partial^2_k\phi(x) \\
&+  \frac{1}{c}\partial_k g(x/c)\exp\left(\frac{\phi(x)}{2}\right)\partial_k \phi(x) + \frac{1}{2}g(x/c)\exp\left(\frac{\phi(x)}{2}\right)(\partial_k\phi(x))^2 \\
&= \sum_{k=1}^d  -\frac{1}{c^2}\partial^2_k g(x/c)\exp\left(\frac{\phi(x)}{2}\right) + \frac{1}{4}g(x/c)\exp\left(\frac{\phi(x)}{2}\right)(\partial_k\phi(x))^2 - \frac{1}{2}g(x/c)\exp\left(\frac{\phi(x)}{2}\right)\partial^2_k\phi(x) \\
&= \left(-\frac{1}{c^2}\Delta g(x/c) +\left(\frac{1}{4}|\nabla\phi(x)|^2 - \frac{1}{2}\Delta\phi(x)\right)g(x/c)\right)\exp\left(\frac{\phi(x)}{2}\right).
\end{align*}
So
\begin{align*}
(Hg)(x) &= -\frac{1}{c^2}\Delta g(x) +\left(\frac{1}{4}|\nabla\phi(cx)|^2 - \frac{1}{2}\Delta\phi(cx)\right)g(x)\\
&= \left(\left(-\frac{1}{c^2}\Delta + C_\phi(\cdot)\right)g(\cdot)\right)(x).
\end{align*}
This can be expressed directly in the standard Weyl calculus as $h(X,D)$, for $h(x, \xi) := \frac{\xi^2}{c^2} + C_\phi (x)$ with the same interpretation as in the statement of this theorem. But $h(Q,P)f = U_2^{-1} h(X,D) U_2f = Lf$, and so we are done.
\end{proof}

\section{$L^p$ bounds on Weyl Pseudodifferential Operators}

In this section we investigate properties of the generalised Weyl calculus associated with a potential $\phi$, insofar as they relate to the functional calculus for the corresponding Ornstein-Uhlenbeck operator. In a later section we will return to the study of the symbol calculus.

Our symbols will correspond to certain functions $a\in L^\infty(\mathbb{R}^{2d})$, which shall be denoted $a(x, \xi)$ for $x, \xi \in \mathbb{R}^d$. To define our symbol class, we need a few intermediate definitions.

\begin{defn}\label{defnD}
Define the Banach space of $L^1$ dominated functions $\mathcal{D}_{\infty,1} = \{ g \in L^\infty(\mathbb{R}^d; L^1(\mathbb{R}^d)); \exists G \in L^1(\mathbb{R}^d), |g(x,k)|<G(k) \text{ for a.e. } x,k \in \mathbb{R}^d\}$, equipped with the norm
\[||g||_{\mathcal{D}_{\infty,1}} = \inf\left\{ ||G||_{L^1(\mathbb{R}^d)}; G \in L^1(\mathbb{R}^d), |g(x,k)|<G(k) \text{ for a.e. } x,k \in \mathbb{R}^d\right\}.\]
\end{defn}

We won't prove that $\mathcal{D}_{\infty,1}$ is a Banach space, although it is easy. Subadditivity and homogeneity of $||\cdot||_{\mathcal{D}_{\infty,1}}$ is obvious. That $||\cdot||_{\mathcal{D}_{\infty,1}}$ is positive definite can be seen by noting it is bounded below by the $L^\infty(\mathbb{R}^d; L^1(\mathbb{R}^d))$ norm. Checking that $\mathcal{D}_{\infty,1}$ is complete is a standard exercise in telescoping sums, and using the fact that $L^\infty(\mathbb{R}^d; L^1(\mathbb{R}^d))$ is complete.

\begin{defn}
For $a\in L^\infty(\mathbb{R}^{2d})$, define $I_a:\mathbb{R}^d \to \mathcal{S}'(\mathbb{R}^d)$, via the action at $x \in \mathbb{R}^d$ and for $\varphi \in \mathcal{S}(\mathbb{R}^d)$ as
\[\left\langle I_a(x), \varphi \right\rangle = \int_{\mathbb{R}^d} a(x, \xi)\varphi(\xi)d\xi.\]
\end{defn}

We will make use of the Fourier transform $\mathcal{F}$ acting on tempered distributions $\sigma \in \mathcal{S}'(\mathbb{R}^d)$, which we normalise such that $\langle \mathcal{F}\sigma, \varphi \rangle = \langle \sigma, \mathcal{F}^*\varphi \rangle$ for all $\varphi \in \mathcal{S}(\mathbb{R}^d)$, where
\[\mathcal{F}^*\varphi(k) = (2\pi)^{-\frac{d}{2}}\int_{\mathbb{R}^d} \varphi(\xi)\exp(i\xi k)d\xi.\]

Now we can define our symbol class.

\begin{defn}\label{defnHS0M}
Fix $M:\mathbb{R}^d \to \mathbb{R}^d$. The space $HS_0(M)$ (standing for Holomorphic Strip) is a subspace of $L^\infty(\mathbb{R}^{2d})$, with
\[HS_0(M) = \left\{ 
\begin{array}{c}
 a \in L^\infty(\mathbb{R}^{2d}); \exists ! g_a \in \mathcal{D}_{\infty,1}, \forall x \in \mathbb{R}^d, \forall \varphi \in \mathcal{S}(\mathbb{R}^d),\\
   \left\langle \mathcal{F}(I_a(x)), \varphi \right\rangle = (2\pi)^\frac{d}{2}\int_{\mathbb{R}^d}\exp\left(-|M(x)k|\right)g_a(x,k)\varphi(k)dk
\end{array} \right\}\]
and norm defined by
\[||a||_{HS_0(M)} = ||g_a||_{\mathcal{D}_{\infty,1}}.\]
For $a \in HS_0(M)$, we define $\mathcal{F}_2a$ to be the measurable function $\mathbb{R}^{2d} \to \mathbb{C}$ with action 
\[(x, k) \mapsto (2\pi)^\frac{d}{2} \exp\left(-|M(x)k|\right)g_a(x,k).\]
\end{defn}

In some sense, $\mathcal{F}_2a(x,k)$ is the ``Fourier transform in the second ($\xi$) variable" with $k$ the variable dual to $\xi$, which is our reason for using this notation. It is easy to verify that if $a(x, \cdot) \in L^1(\mathbb{R}^d)$ for each $x\in \mathbb{R}^d$, then $\mathcal{F}_2a$ is indeed the Fourier transform in the second variable for each fixed $x$. However, not all symbols we will consider will be integrable in the second variable, which makes our definition via the space of tempered distributions useful.

We will also often refer to the $g_a \in \mathcal{D}_{\infty,1}, G_a \in L^1(\mathbb{R}^d)$ associated with $a \in HS_0(M)$, by which we mean the unique such $g_a$ as seen in Definition \ref{defnHS0M}, and $G_a$ dominating $g_a$ as in Definition \ref{defnD}. The main reason for using such a complicated definition is that the integral kernel of $a(Q,P)$ for $a \in \mathbb{S}(\mathbb{R}^{2d})$ from Theorem \ref{thmKernel} is closely related to the Fourier transform in the second variable of the Schwartz function $a$, as in the following lemma.

\begin{lemma}\label{lemmaKernel}
Fix a potential $\phi \in C^2$, and let $(Q, P)$ be the associated generalised Weyl pair. For $a \in \mathcal{S}(\mathbb{R}^{2d})$, the integral kernel of the operator $a(Q,P)$ is given by
\[k(y,x) = \frac{1}{(2\pi)^{d/2}c^d} \mathcal{F}_2 a\left(\frac{x+y}{2c}, \frac{x-y}{c}\right) \exp\left(\frac{1}{2}(\phi(y)-\phi(x))\right).\]
\end{lemma}

\begin{proof}
From Theorem \ref{thmKernel}, we have for $f \in \mathcal{C}_\phi$, $y \in \mathbb{R}^d$,
\begin{align*}
(a(Q,P)f)(y) &= \frac{1}{(2\pi)^{d}c^d} \int_{\mathbb{R}^{2d}}a\left(\frac{x+y}{2c}, \xi\right) \exp\left(-i\xi \left(\frac{x-y}{c}\right)\right) \exp\left(\frac{1}{2}(\phi(y)+\phi(x))\right)f(x) d\xi d\mu(x) \\
&= \frac{1}{(2\pi)^{d}c^d} \int_{\mathbb{R}^{2d}}a\left(\frac{x+y}{2c}, \xi\right) \exp\left(-i\xi \left(\frac{x-y}{c}\right)\right) \exp\left(\frac{1}{2}(\phi(y)-\phi(x))\right)f(x) d\xi dx \\
&= \frac{1}{(2\pi)^{d/2}c^d} \int_{\mathbb{R}^{d}}\mathcal{F}_2a\left(\frac{x+y}{2c}, \frac{x-y}{c}\right) \exp\left(\frac{1}{2}(\phi(y)-\phi(x))\right)f(x) dx, \\
\end{align*}
noting that the order of integration can be changed, and the integration in $\xi$ carried out, as everything converges absolutely.
\end{proof}

We make use of Lemma \ref{lemmaKernel} to extend the generalised Weyl calculus to $HS_0(M)$, as in the following definition:

\begin{defn}\label{defnWeylHS0M}
For any function $M:\mathbb{R}^d \to \mathbb{R}^d$, and any potential $\phi \in C^2$, define for $a \in HS_0(M)$ the operator $a(Q,P):\mathcal{C}_\phi \to L^1_{loc}(\lambda)$ via the action
\[(a(Q,P)f)(y) = \frac{1}{(2\pi)^{d/2}c^d} \int_{\mathbb{R}^d} \mathcal{F}_2 a\left(\frac{x+y}{2c}, \frac{x-y}{c}\right) \exp\left(\frac{1}{2}(\phi(y)-\phi(x))\right)f(x)dx.\]
\end{defn}

To ensure this definition makes sense, we check that for $f \in \mathcal{C}_\phi$, $a(Q,P)f \in L^1_{loc}(\mathbb{R}^d)$.

\begin{proof}
Using the fact that $a \in HS_0(M)$ for some $M:\mathbb{R}^d \to \mathbb{R}^d$, there exists a $G \in L^1(\mathbb{R}^d)$ such that $\left|\mathcal{F}_2a(x,k)\right| \leq G(k)$ for a.e. $x,k \in \mathbb{R}^d$. Hence we find for a.e. $y \in \mathbb{R}^d$
\begin{align*}
\left|(a(Q,P)f)(y)\right| &= \frac{1}{(2\pi)^{d/2}c^d} \left|\int_{\mathbb{R}^{d}}\mathcal{F}_2a\left(\frac{x+y}{2c}, \frac{x-y}{c}\right) \exp\left(\frac{1}{2}(\phi(y)-\phi(x))\right)f(x) dx\right| \\
&\leq \frac{1}{(2\pi)^{d/2}c^d} \int_{\text{supp}(f)}G\left(\frac{x-y}{c}\right) \exp\left(\frac{1}{2}(\phi(y)-\phi(x))\right)|f(x)| dx
\end{align*}
Since $f \in \mathcal{C}_\phi$, $\text{supp}(f)$ is compact and $||f||_{L^\infty(\lambda)}$ is bounded. Since $\phi$ is continuous, it will be bounded on $\text{supp}(f)$, by $C>0$ say. So
\begin{align*}
\left|(a(Q,P)f)(y)\right| &\leq \frac{1}{(2\pi)^{d/2}c^d} \exp(C/2)||f||_{L^\infty(\lambda)}\exp\left(\phi(y)/2\right)\int_{\text{supp}(f)}G\left(\frac{x-y}{c}\right)  dx \\
&\leq \frac{1}{(2\pi)^{d/2}c^d} \exp(C/2)||f||_{L^\infty(\lambda)}||G||_{L^1(\lambda)}\exp\left(\phi(y)/2\right).\\
\end{align*}
Since $\phi$ is continuous, $\left(y \mapsto \exp\left(\phi(y)/2\right) \right) \in L^1_{loc}(\lambda)$, and thus so is $a(Q,P)f$.
\end{proof}

We now show our main theorem: that for $a \in HS(M)$ for the correct $M$, $a(Q,P)$ extends to a bounded operator on $L^p(\mu)$. The correct $M$ is as follows.

\begin{defn}
A pair $(M,\epsilon)$  consisting of a measurable function $M:\mathbb{R}^d \to \mathbb{R}^d$ and a number $\epsilon \geq 0$ is a called a valid growth pair for $\phi \in C^2(\mathbb{R}^d)$ and $p \in [1, \infty]$ if for all $x,y \in \mathbb{R}^d$,
\[\left|\left( \left|\frac{1}{2}-\frac{1}{p}\right|\left|\phi(x)-\phi(y)\right| - \left|\left(\frac{x-y}{c}\right) M\left(\frac{x+y}{2c}\right)\right|\right)\right| \leq \epsilon.\]
\end{defn}

\begin{thm}\label{thmBoundedness0}
Fix a potential $\phi \in C^2$, let $(Q, P)$ be the associated generalised Weyl pair, and fix $p \in [1,\infty]$. Suppose there exists a valid growth pair $(M,\epsilon)$ for $\phi$ and $p$. Then for $a \in HS_0(M)$ the operator $a(Q,P)$, defined as in Definition \ref{defnWeylHS0M}, extends to a bounded operator on $L^p(\mu)$ and
\[||a(Q,P)||_{B(L^p(\mu))} \leq e^\epsilon||a||_{HS_0(M)}.\]
That is, the generalised Weyl calculus extends to a bounded linear map $HS_0(M) \to B(L^p(\mu))$.
\end{thm}

\begin{proof}
Let $U_p:L^p(\mu) \to L^p(\lambda)$ be the isometry from Definition \ref{defnUp}. Then $a(Q,P)$ has a bounded extension to $L^p(\mu)$ if and only if $U_pa(Q,P)U_p^{-1}$ has a bounded extension to $L^p(\lambda)$, in which case $||a(Q,P)||_{B(L^p(\mu))} = ||U_pa(Q,P)U_p^{-1}||_{B(L^p(\lambda))}$. We will use a Young's convolution inequality argument to show that under our conditions $U_pa(Q,P)U_p^{-1}$ extends boundedly to $L^q(\lambda)$ for all $q \in [1, \infty]$ with norm at most $e^\epsilon ||a||_{HS_0(M)}$. We then remove the isometries on $L^p$ to obtain the desired result.

From Definition \ref{defnWeylHS0M}, $U_pa(Q,P)U_p^{-1}$ can be expressed as an integral operator on $L^q(\lambda)$ with kernel $k:\mathbb{R}^{2d} \to \mathbb{C}$ where for $x, y \in \mathbb{R}^d$,
\[k(y,x) = \frac{1}{(2\pi)^{d/2} c^d} \mathcal{F}_2a\left(\frac{x+y}{2c}, \frac{x-y}{c}\right) \exp\left(\left(\frac{1}{2}-\frac{1}{p}\right)(\phi(y)-\phi(x))\right).\]
Using the definition of $\mathcal{F}_2 a$ and $HS_0(M)$, we find
\begin{align*}
k(y,x) &= \frac{1}{(2\pi)^{d/2} c^d} \mathcal{F}_2a\left(\frac{x+y}{2c}, \frac{x-y}{c}\right)  \exp\left(\left(\frac{1}{2}-\frac{1}{p}\right)(\phi(y)-\phi(x))\right)\\
&= \frac{1}{c^d}g_a\left(\frac{x+y}{2c},\frac{x-y}{c}\right)\exp\left(\left(\frac{1}{2}-\frac{1}{p}\right)(\phi(y)-\phi(x))-\left|\left(\frac{x-y}{c}\right) M\left(\frac{x+y}{2c}\right)\right|\right),\\
\end{align*}
for a unique $g_a \in \mathcal{D}_{\infty,1}$.
Our assumption implies $\left(\frac{1}{2} -\frac{1}{p}\right)(\phi(x)-\phi(y)) - \left|\left(\frac{x-y}{c}\right) M\left(\frac{x+y}{2c}\right)\right| \leq \epsilon$, which we incorporate to find
\begin{align*}
|k(y,x)| &\leq   e^\epsilon \frac{1}{c^d}\left|g_a\left(\frac{x+y}{2c},\frac{x-y}{c}\right)\right|\\
&\leq e^\epsilon \frac{1}{c^d}G_a\left(\frac{x-y}{c}\right),
\end{align*}
for some $G_a \in L^1(\mathbb{R}^d)$. Young's convolution inequality implies extendability and that
\begin{align*}
||U_pa(Q,P)U_p^{-1}||_{B(L^q(\lambda))} &\leq e^\epsilon\frac{1}{c^d}\int_{\mathbb{R}^d}G_a\left(\frac{x}{c}\right)dx \\
&= e^\epsilon\int_{\mathbb{R}^d}G_a\left(x\right)dx.
\end{align*}
Taking infimum over all $G_a$ dominating $g_a$ in the sense of Definition \ref{defnD} gives $||U_pa(Q,P)U_p^{-1}||_{B(L^q(\lambda))} \leq e^\epsilon||a||_{HS_0(M)}$. Removing the isometries on $L^p$ we obtain our desired result.
\end{proof}

\begin{rem}[Existence of a Valid Growth Pair]\label{remAffine}
If $p=2$, the function taking value $0$, and $\epsilon=0$ is a valid growth pair for any $\phi \in C^2(\mathbb{R}^d)$.
If $\phi(x) = xN(x) + lx + \epsilon(x)$, where $N$ is a linear map $\mathbb{R}^d \to \mathbb{R}^d$, $l \in \mathbb{R}^d$, and $\epsilon$ is a bounded $C^2$ function, then $M$ can be taken to be an affine function of $x$ with real-symmetric linear part, depending only on $c, p, N$ and $l$. We will not make such an assumption until Section \ref{secMoreHSM}, and so we will keep $M$ as a general measurable function from $\mathbb{R}^d$ to itself unless otherwise specified.

Note that for $\phi$ of the form $x \mapsto xN(x) + lx + \epsilon(x)$ with $\epsilon(x)$ a bounded $C^2$ function, an operator will be bounded on $L^p(\mu)$ if and only if it is bounded on $L^p(\tilde{\mu})$, where $\tilde{\mu}$ is associated with $\tilde{\phi}(x) = xN(x) + lx$, as the two measures are equivalent and the Radon-Nikodym derivative of one with respect to the other is $\exp(\pm \epsilon(x))$, which is bounded above and below by positive constants. Thus boundedness of operators in our generalised Weyl calculus should only depend on $c, p, N$ and $l$. However, the operators which we should care about (such as the relevant Ornstein-Uhlenbeck operator $L$), will depend on all of $\phi$, not just its unbounded terms.
\end{rem}

\begin{rem}[Extension]
It should be noted that the way this functional calculus is defined is quite different to the standard methods. Rather than a convergence lemma or a density argument, we have found an integral operator expression for ``nice" symbols, and then extended to a large class of symbols for which the integral representation can be made sense of. Thus when we say that the generalised Weyl calculus extends to $HS_0(M)$, we mean both that for $a \in HS_0(M)$, $a(Q,P):\mathcal{C}_\phi \to L^1_{loc}(\lambda)$ as defined in Definition \ref{defnWeylHS0M} extends uniquely to a bounded operator on $L^p(\mu)$ (by density of $\mathcal{C}_\phi$ as in Theorem \ref{thmCore}), and also that if $a \in \mathcal{S}(\mathbb{R}^{2d}) \cap HS_0(M)$, then the expressions for $a(Q,P)$ from Definitions \ref{defnWeylCalc} and \ref{defnWeylHS0M} agree (by Lemma \ref{lemmaKernel}).
\end{rem}

\begin{rem}\label{remInitial}
Suppose that $a:\mathbb{R}^{2d} \to \mathbb{C}$ is such that for each $x \in \mathbb{R}^d$, $a(x, \xi + i \eta)$ can be extended to a holomorphic function for $\eta$ in $B(0, |M(x)|)$, and such that there exists a constant $K>0$ such that for any multiindex $\alpha$ with $|\alpha| \leq (d+1)$ we have
\[\sup_{u \in \mathbb{R}^d, \eta \in B(0, |M(u)|)} \int_{\mathbb{R}^{d}}\left|\partial_\xi^{(\alpha)}a\left(u, \xi+i\eta\right)\right| d\xi \leq K.\]
Then $a \in HS_0(M)$.
However, this integrability condition is much stronger than what we actually require for admission into $HS_0(M)$. Although, as we shall see, some ``pseudo-holomorphic" nature is always apparent for symbols in $HS_0(M)$.
\end{rem}

\begin{rem}[Holomorphic Nature of $HS_0(M)$]\label{remHolomorphicExtension}

It is well-known that if a function $f \in L^\infty(\mathbb{R}^d)$ has Fourier transform of the form
\[(\mathcal{F}f)(\xi) = (2\pi)^{\frac{d}{2}}\exp \left(-a|\xi|\right) g(\xi)\]
for some $a>0$, $g \in L^1(\mathbb{R}^d)$, then $f$ almost everywhere agrees with the restriction to $\mathbb{R}^d$ of a function holomorphic on the cylinder $\{\xi + i \eta; |\eta| < a\}$ (this can be verified via a change of contour in the integral expression for the inverse Fourier transform). Furthermore, this holomorphic extension has a continuous extension to the closure of the cylinder, and the supremum norm of said continuous extension is bounded by $||g||_{L^1(\mathbb{R}^d)}$.

There is an analogous statement for elements of $HS_0(M)$. Fix $a \in HS_0(M)$, $x, \xi \in \mathbb{R}^d$ and $t \in [-1, 1]$. Then we find
\begin{align*}
\frac{1}{(2\pi)^\frac{d}{2}} \left|\int_{\mathbb{R}^d} \mathcal{F}_2a(x, k)\exp(i k (\xi + i t M(x)) ) dk\right| &\leq \frac{1}{(2\pi)^\frac{d}{2}}\int_{\mathbb{R}^d} |\mathcal{F}_2a(x, k)|\exp(- t k M(x)) ) dk \\
&= \int_{\mathbb{R}^d} \exp\left(- |k M(x)| - t k M(x) \right)|g_a(x,k)|dk \\
&\leq \int_{\mathbb{R}^d} G_a(k) dk \\
&\leq ||g_a||_{\mathcal{D}_{\infty,1}} \\
&= ||a||_{HS_0(M)},
\end{align*}
where the infimum is over all $G_a$ dominating $g_a$ as in Definition \ref{defnD}. By the Fourier inversion formula, the above agrees $\xi$-almost everywhere with $a(x,\xi)$ when $t=0$. We thus take the above as a definition of an extension of $a$ to the set 
\[D_M = \{ (x, \xi+i\eta) \in \mathbb{R}^d \times \mathbb{C}^d; \xi \in \mathbb{R}^d, \exists t \in [-1,1] \text{ s.t. } \eta = tM(x).\]
We denote this extension at $x, \xi, \eta$ as $a(x, \xi+ i\eta)$. For each fixed $x$, continuity as a function of $\xi$ and $t$ follows by the dominated convergence theorem. If $d=1$, for fixed $x$ this extension is holomorphic as a function of $\xi + i\eta \in \mathbb{C}$, which also follows from the DCT.

If $d \geq 2$ and fixed $x\in \mathbb{R}^d$, it does not make sense to speak of holomorphy of $a(x, \cdot)$ due to its domain not being an open subset of $\mathbb{C}^d$. However, for fixed $x \in \mathbb{R}^d$, $(\xi, t) \mapsto a(x, \xi +itM(x))$ will be a real-analytic function on $\mathbb{R}^d \times (-1,1)$ and will satisfy some modified form of the Cauchy-Riemann equations.

That functions $a \in HS_0(M)$ possess for each fixed $x\in \mathbb{R}^d$ such a ``pseudo-holomorphic" extension to the strip $\{\xi \in \mathbb{R}^d, \exists t \in [-1,1] \text{ s.t. } \eta = tM(x)\}$ is where the name Holomorphic Strip originated.
\end{rem}

\begin{rem}[Comparison to the standard symbol classes]
We should compare the symbol class $HS_0(M)$ to standard symbol classes giving rise to bounded operators through the Weyl calculus. If we take $\phi(x) = 0$ our space of functions is $L^p(\lambda)$ and our generalised Weyl pair is the standard one, in which there are many known symbol classes giving rise to bounded operators (see for instance \cite{Stein}, Chapter 6). These classes typically assume boundedness and decay in the $\xi$ variable of sufficiently many derivatives of $a(x,\xi)$, and allow for some singular integral operators. If $\phi(x)=0$ we can take $M=0$, in which case $HS_0(0)$ will be the space of functions whose Fourier transform in the second variable is dominated by an integrable function, thus not including singular integral operators. This implies boundedness, continuity in $\xi$ and decay of $a(x,\xi)$ (by the Riemann-Lebesgue Lemma), but does not give a rate of decay or any differentiability. Similarly, in the case $p=2$ but $\phi(x) \neq 0$, our generalised Weyl calculus is unitarily equivalent to the standard Weyl calculus, and we can again take $M=0$.

Alternatively, when $\phi(x) \neq 0$ and $p \neq 2$, we find the relevant $M$ is non-zero and so by Remark \ref{remHolomorphicExtension}, symbols in $HS_0(M)$ must have  pseudo-holomorphic extendability condition. In this case, there is no isometry back to $L^p(\lambda)$, mapping the associated Weyl pair to the standard Weyl pair. At first sight, this seems infinitely worse than the standard symbol classes. However, this may be the best that can be done. The classical Ornstein-Uhlenbeck operator, associated with $\phi(x)=\frac{x^2}{2}$, is known to only have bounded functional calculus which is holomorphic (see \cite{GMMST}), and so if we expect to be able to study the functional calculus of the classical Ornstein-Uhlenbeck operator via the associated generalised Weyl calculus we should be forced into accepting some sort of holomorphic extendability condition on the symbols which give rise to bounded operators.

For the specific case $\phi(x)=\frac{x^2}{2}$, we can factorise the exponential term $\exp\left(\left(\frac{1}{2}-\frac{1}{p}\right)(\phi(y)-\phi(x))\right)$ of the integral kernel of $U_p\circ a(Q,P) \circ U_p^{-1}$ into a function of $\frac{x+y}{2}$ and $x-y$. Using this, for any $U_p\circ a(Q,P) \circ U_p^{-1}$ for $(Q,P)$ associated with $\phi(x)$, we can find a symbol $\tilde{a}$ such that  $U_p\circ a(Q,P) \circ U_p^{-1} = \tilde{a}(X,D)$, where $(X,D)$ is the standard Weyl pair. Thus we could derive boundedness of $a(Q,P)$ by checking when $\tilde{a}$ satisfies the standard symbol estimates of classical pseudodifferential operator theory. However, by carrying out this calculation formally, we find that for $x, \xi \in \mathbb{R}^d$, $\tilde{a}(x,\xi) = a\left(\frac{x}{c}, c\xi + i c \left(\frac{1}{2}-\frac{1}{p}\right)x\right)$, defined as in Remark \ref{remHolomorphicExtension}. So for such an argument to work, we would need standard symbol estimates on the ``boundary" of this extension. It is clear that this method would still lead to some strong restrictions on symbols.
\end{rem}

We can push the techniques used to prove Theorem \ref{thmBoundedness0} ever so slightly to prove the following R-boundedness theorem. See \cite{AnalysisInBanachSpacesV2} for the theory of R-boundedness.

\begin{thm}\label{thmRBoundedness0}
Fix $p \in (1, \infty)$, a potential $\phi$, and suppose there exists a valid growth pair $(M,\epsilon)$ for $\phi$ and $p$.
Let $A \subset HS_0(M)$, and $G \subset \mathcal{D}_{\infty,1}$ be the set of $g_a$ corresponding to each $a \in A$ as in Definition \ref{defnHS0M}. Suppose that for each $g_a \in G$ we can choose a dominating $G_a \in L^1(\mathbb{R}^d)$ such that the supremum over our selections of the quantity
\[\int_{\mathbb{R}^d} \esssup_{|y| \geq |x|} |G_a(y)| dx\]
is finite. Then $A(Q,P) = \{a(Q,P), a \in A\}$ is R-bounded on $L^p(\mu)$.
\end{thm}

\begin{proof}
We apply the same technique as was used in Theorem \ref{thmBoundedness0}, introducing the isometries $U_p: L^p(\mu) \to L^p(\lambda)$ from Definition \ref{defnUp}. We will show that the set $U_p \circ A(Q,P) \circ U_p^{-1}$ is R-bounded on $L^q(\lambda)$ for all $q \in (1, \infty)$, specifically for $q=p$, in which case we obtain R-boundedness of $A(Q,P)$ on $L^p(\mu)$ by removing the isometries.
As in Theorem \ref{thmBoundedness0}, the integral kernel of an operator $U_p \circ a(Q,P) \circ U_p^{-1}$ is given by
\[ k_a(y,x) = \frac{1}{(2\pi)^{d/2} c^d} \mathcal{F}_2a\left(\frac{x+y}{2c}, \frac{x-y}{c}\right) \exp\left(\left(\frac{1}{2}-\frac{1}{p}\right)(\phi(y)-\phi(x))\right) \]
and so for all $x, y \in \mathbb{R}^d$,
\begin{align*}
|k_a(y,x)| &\leq \frac{1}{c^d} \left|g_a\left(\frac{x+y}{2c},\frac{x-y}{c}\right)\right| \exp\left(-\left|\left(\frac{x-y}{c}\right)M\left(\frac{x+y}{2c}\right)\right|\right) \exp\left(\left(\frac{1}{2}-\frac{1}{p}\right)(\phi(y)-\phi(x))\right) \\
&\leq e^\epsilon \frac{1}{c^d} G_a\left(\frac{x-y}{c}\right)
\end{align*}
So each $U_p \circ a(Q,P) \circ U^{-1}_p$ has kernel dominated by a convolution, namely convolution against $e^\epsilon \frac{1}{c^d} G_a\left(\frac{\cdot}{c}\right)$. Hence R-boundedness of the set of convolution operators $\{f \mapsto e^\epsilon\frac{1}{c^d} G_a\left(\frac{\cdot}{c}\right) * f\}$ on $L^q(\mathbb{R}^d)$ will imply R-boundedness of $U_p \circ A(Q,P) \circ U^{-1}_p$ on $L^q(\mathbb{R}^d)$ (see Proposition 8.1.10 of \cite{AnalysisInBanachSpacesV2}, and note that in the proof of said proposition the fixed positive operator can be replaced by an R-bounded family of positive operators). Applying Proposition 8.2.3 of \cite{AnalysisInBanachSpacesV2} and our assumptions shows that $\{f \mapsto e^\epsilon\frac{1}{c^d} G_a\left(\frac{\cdot}{c}\right) * f\}$ is R-bounded on $L^q(\lambda)$ for all $q \in (1,\infty)$. Therefore, $A(Q,P)$ is R-bounded on $L^p(\mu)$.
\end{proof}

\section{An Application}\label{secHInfty}

In this section we will use our generalised Weyl calculus developed in the previous section to show that for $\phi$ of the form $\phi(x) = x N(x) + lx$, where $N:\mathbb{R}^d \to \mathbb{R}^d$ is a positive-definite real-symmetric linear operator and $l \in \mathbb{R}^d$ is a vector, the associated Ornstein-Uhlenbeck operator has bounded $H^\infty(\Sigma_{\theta_p})$ functional calculus on $L^p(\mu)$, where $\sin (\theta_p) = M_p := |1-\frac{2}{p}|$. This result generalises that of \cite{GMMST}, which shows that the given angle is optimal for the bounded $H^\infty$ functional calculus of the classical Ornstein-Uhlenbeck operator (corresponding to $N$ as half the identity and $l=0$). In this proof, the use of the generalised Weyl calculus can be a posteriori removed, leading to a strikingly simple proof in the classical case (see \cite{HARRIS2019}). We include the argument here to show that Theorem \ref{thmBoundedness0} has important consequences despite the simplicity of its proof and the complexity of the definition of $HS_0(M)$.

\begin{thm}\label{corHinftyGoodAngle}
For $p \in (1, \infty)$, $\phi(x) = x N(x) + lx$ where $N:\mathbb{R}^d \to \mathbb{R}^d$ is a positive definite real-symmetric linear operator and vector $l \in \mathbb{R}^d$, the associated Ornstein-Uhlenbeck operator has bounded $H^\infty(\Sigma_{\theta_p})$ functional calculus on $L^p(\mu)$, where $\sin (\theta_p) = M_p := |1-\frac{2}{p}|$.
\end{thm}

Our proof is based on the well-known result that in $L^p$ spaces, if an operator is known to have a bounded $H^\infty$ functional calculus of some angle, the optimal angle of the $H^\infty$ functional calculus of the operator is equal to its optimal angle of R-sectoriality (see \cite{AnalysisInBanachSpacesV2} for the theory of R-sectoriality, and its Theorem 10.7.13 for a proof of the stated result). We have already seen that any of the Ornstein-Uhlenbeck operators considered in this paper automatically have $H^\infty$ functional calculus of some angle (see Corollary \ref{corSomeHinfty}), so all that we need to do is optimise the angle. Our proof that the angle of R-sectoriality of the relevant Ornstein-Uhlenbeck operator is equal to $\theta_p$ uses Theorem 10.3.3 of \cite{AnalysisInBanachSpacesV2}, which states an equivalence between an operator $A$ being R-sectorial of angle $\theta<\frac{\pi}{2}$ and $-A$ being the generator of an analytic semigroup of angle $\frac{\pi}{2}-\theta$ which is R-bounded on each smaller sector. We prove the required analytic extendability and R-boundedness of the relevant Ornstein-Uhlenbeck semigroup by using the following generalisation of Theorem \ref{TheoremSGEasy} to transfer to a study of the generalised Weyl calculus, and Theorem \ref{thmRBoundedness0} to deduce the required R-boundedness result.

\begin{thm}\label{TheoremSGHard}
Suppose $\phi(x) = xN(x) + lx$ for a positive semi-definite real-symmetric linear map $N: \mathbb{R}^d \to \mathbb{R}^d$ and vector $l \in \mathbb{R}^d$. Then for the associated Ornstein-Uhlenbeck operator $L$ we have for all $t >0$ (initially defined as a map $\mathcal{C}_\phi \to L^1_{loc}(\mathbb{R}^d)$)
\[\exp(-tL) = a_t^{N,l}(Q,P)\]
where
\begin{equation*}
a_t^{N,l}(x, \xi) = \det\left(\cosh(tN)^{-1}\exp(tN)\right)\exp\left(-\frac{1}{c^2}\xi N_t(\xi) -\left(cN(x) + \frac{l}{2}\right)N_t\left(cN(x) + \frac{l}{2}\right)\right)
\end{equation*}
and
\[N_t = N^{-1}\tanh(tN),\]
with all functions of $N$ interpreted via the functional calculus of a real-symmetric operator on $\mathbb{R}^d$ with the standard inner product (if $0 \in \sigma(N)$, we set $N_t$ to act as multiplication by $t$ on the $0$-eigenspace, which can be motivated by noting that for fixed $t >0$, $n \mapsto \frac{\tanh{tn}}{n}$ has a unique entire analytic extension, with value $t$ at $n=0$).
\end{thm}

\begin{proof}
Note that $L$ with domain $\mathcal{C}_\phi$ is unitarily equivalent to $-\frac{1}{c^2}\Delta + C_\phi(x)$ with domain $C_c^\infty(\mathbb{R}^d)$, where $C_\phi(x)=\frac{1}{4}|\nabla\phi(cx)|^2 - \frac{1}{2}\Delta \phi(cx)$, as is seen in the proof of Theorem \ref{thmH}. Note that as $\phi(x)$ is a second degree polynomial in the components of $x$, $C_\phi(x)$ will also be a second degree polynomial in the components of $x$. The work of H\"{o}rmander in \cite{UnpublishedHormander} gives an explicit representation for the classical Weyl symbol for the semigroup generated by operators of the form $-\Delta+V(x)$ where $V$ is a degree two polynomial in the components of $x$. Using H\"{o}rmander's formula in our particular case and noting the joint unitary equivalence of $L$ with $-\frac{1}{c^2}\Delta + C_\phi(x)$ and $a(Q,P)$ with $a(X,D)$, we obtain the claimed expression for $a_t^{N,l}$.
\end{proof}

To deduce the desired $H^\infty$ functional calculus result, we need only show that the relevant Ornstein-Uhlenbeck semigroup has an analytic extension to a sector of the correct angle, and that it is R-bounded on each smaller sector. We will in fact show a lot more with no more effort. The function $t \mapsto N_t$ is analytic and can clearly be extended to $\mathbb{C}\backslash \bigcup_{n \in \sigma(N)} \frac{i\pi}{n}(2\mathbb{Z}+1)$. We will consider the analytic extension $z \mapsto N_z$ on domains of the form
\begin{equation}\label{eqnEN}
E^N := \{z \in \mathbb{C}; \sigma(N_z) \subset \Sigma_{\frac{\pi}{2}-\theta_p}\}\backslash \bigcup_{n \in \sigma(N)} \frac{i\pi}{n}\left(\mathbb{Z}+\frac{1}{2}\right)
\end{equation}
where $\sin (\theta_p) = M_p := \left|1-\frac{2}{p}\right|$. We will show the Ornstein-Uhlenbeck semigroup extends to an analytic semigroup on the domain $E$. Moreover, we will show that the Ornstein-Uhlenbeck semigroup is R-bounded on sets of the form
\begin{equation}\label{eqnENdash}
E^N_{\epsilon, \delta} := \left\{z \in \mathbb{C}; \sigma(\cos^2(\arg(N_z))) \in (M_p^2+\epsilon, \infty), \text{dist}\left(z, \left(\bigcup_{n \in \sigma(N)} \frac{i\pi}{2n}\mathbb{Z}\right)\backslash \{0\}\right) > \delta\right\}
\end{equation}
for all $\epsilon, \delta>0$. The condition on the spectrum of the cosine of the argument of $N_z$ is a rephrasing of the spectrum of $N_z$ being contained in a sector slightly smaller than $\Sigma_{\theta_p}$, which is a useful form for the proof to come. Note the condition on the distance to $\left(\bigcup_{n \in \sigma(N)} \frac{i\pi}{2n}\mathbb{Z}\right)\backslash \{0\}$ ensures we remain uniformly away from the poles and zeroes of $N_z$, besides $z=0$. We claim that $\Sigma_{\frac{\pi}{2}-\theta_p} \subset E^N$ for any $N$, and that for all $\epsilon'>0$ there exists $\epsilon, \delta>0$ such that $\Sigma_{\frac{\pi}{2}-\theta_p-\epsilon'} \subset E^N_{\epsilon, \delta}$ (see \cite{vNPWeylGaussian} for details of this calculation in the case $N$ is a multiple of the identity, and note that the general case follows by taking intersections over the eigenvalues of $N$). These results combined will imply that the maximal domain of analyticity of the Ornstein-Uhlenbeck semigroup contains the sector $\Sigma_{\frac{\pi}{2}-\theta_p}$, and that it is R-bounded on each smaller sector.

\begin{thm}\label{thmRSectorialGoodAngle}
Suppose $\phi(x) = xN(x) + lx$ for a positive semi-definite real-symmetric linear map $N: \mathbb{R}^d \to \mathbb{R}^d$ and vector $l \in \mathbb{R}^d$. For $p \in (1, \infty)$, the associated Ornstein-Uhlenbeck semigroup has an analytic extension on $L^p(\mu)$ to the domain $E^N$. Furthermore, if $N$ is positive definite, this extension is R-bounded on each domain $E^N_{\epsilon, \delta}$ for all $\epsilon, \delta>0$.
\end{thm}

\begin{proof}
Let $z \in E^N$. By Theorem \ref{TheoremSGHard}, $\exp(-zL) = a_z(Q,P)$ (dropping $N, l$ from the notation), where
\[a_z(x, \xi) = \det\left(\cosh(zN)^{-1}\exp(zN)\right)\exp\left(-\frac{1}{c^2}\xi N_z(\xi) -\left(cN(x) + \frac{l}{2}\right)N_z\left(cN(x) + \frac{l}{2}\right)\right)\]
and $N_z = N^{-1}\tanh(zN)$.
Computing $\mathcal{F}_2a_z(x, k)$ gives:
\begin{align*}
\mathcal{F}_2 a_z(x, k) &= 2^{-\frac{d}{2}}c^d\det\left(N_z^{-\frac{1}{2}}\cosh(zN)^{-1}\exp(zN)\right)\exp\left(-\frac{c^2}{4} k N^{-1}_z(k) -\left(cN(x) + \frac{l}{2}\right)N_z\left(cN(x) + \frac{l}{2}\right)\right)\\
\end{align*}
For $\phi(x) = xN(x) + lx$ and fixed $p \in (1, \infty)$, we claim $M:\mathbb{R}^d \to \mathbb{R}^d$, $x \mapsto M_p \left(c^2N(x) + \frac{cl}{2}\right)$, and $\epsilon=0$ are a valid growth pair. To see this, we have for all $x,y \in \mathbb{R}^d$ 
\begin{align*}
 \left|\frac{1}{2}-\frac{1}{p}\right|\left|\phi(x)-\phi(y)\right| &= \frac{M_p}{2}\left|xN(x) + lx-yN(y) - ly\right|\\
&= \frac{M_p}{2}\left|(x-y)N(x+y) + l(x-y)\right|\\
&= \left| \left(\frac{x-y}{c}\right)M_p\left(c^2N\left(\frac{x+y}{2c}\right) + \frac{cl}{2}\right)\right|.\\
 \end{align*}
 Using this, we can rewrite $\mathcal{F}_2 a_z(x, k)$ as
 \begin{align*}
\mathcal{F}_2 a_z(x, k) &= 2^{-\frac{d}{2}}c^d\det\left(N_z^{-\frac{1}{2}}\cosh(zN)^{-1}\exp(zN)\right)\exp\left(-\frac{c^2}{4} k N^{-1}_z(k) -\frac{1}{c^2M_p^2} M(x)N_z\left(M(x)\right)\right)\\
\end{align*}
We will show that $a_z$ satisfies the conditions of admission into $HS_0(M)$ for this specific $M$, from which Theorem \ref{thmBoundedness0} gives boundedness on $L^p(\mu)$ of $a_z(Q,P)$. The decomposition of $a_z$ as in Definition \ref{defnHS0M} has
\begin{align*}
g_{a_z}(x, k) = \frac{c^d}{2^d\pi^\frac{d}{2}}\det\left(N_z^{-\frac{1}{2}}\cosh(zN)^{-1}\exp(zN)\right)\exp\left(-\frac{c^2}{4} k N^{-1}_z(k) + |k M(x)| -\frac{1}{c^2M_p^2} M(x)N_z\left(M(x)\right)\right)
\end{align*}
We wish to show $g_{a_z} \in \mathcal{D}_{\infty,1}$, so we must dominate in $x$ by an integrable function in $k$. Letting $\iota$ be the sign of $k M(x)$, we find by completing the square in $M(x)$:
\begin{align*}
|g_{a_z}(x, k)| &= \frac{c^d}{2^d\pi^\frac{d}{2}}\left|\det\left(N_z^{-\frac{1}{2}}\cosh(zN)^{-1}\exp(zN)\right)\right|\\
&\exp\left(-\frac{c^2}{4} k \Re(N^{-1}_z)(k) + \iota k M(x) -\frac{1}{c^2M_p^2} M(x)\Re(N_z)\left(M(x)\right)\right)\\
&= \frac{c^d}{2^d\pi^\frac{d}{2}}\left|\det\left(N_z^{-\frac{1}{2}}\cosh(zN)^{-1}\exp(zN)\right)\right|\\
&\exp\left(-\frac{c^2}{4} k \Re(N^{-1}_z)(k) + \frac{c^2M_p^2}{4} k\Re(N_z)^{-1}(k) - \left(\frac{1}{cM_p} \Re(N_z)^{\frac{1}{2}}(M(x)) - \frac{\iota cM_p}{2}\Re(N_z)^{-\frac{1}{2}}(k) \right)^2\right)\\
&= \frac{c^d}{2^d\pi^\frac{d}{2} }\left|\det\left(N_z^{-\frac{1}{2}}\cosh(zN)^{-1}\exp(zN)\right)\right|\\
&\exp\left(-\frac{c^2}{4} k \left( \Re(N^{-1}_z) - M_p^2 \Re(N_z)^{-1}\right)(k) - \left(\frac{1}{cM_p} \Re(N_z)^{\frac{1}{2}}(M(x)) - \frac{\iota cM_p}{2}\Re(N_z)^{-\frac{1}{2}}(k) \right)^2\right),\\
\end{align*}
where by $\left(\frac{1}{cM_p} \Re(N_z)^{\frac{1}{2}}(M(x)) - \frac{\iota cM_p}{2}\Re(N_z)^{-\frac{1}{2}}(k) \right)^2$ we mean the inner product of the contents of the brackets with itself, which is non-negative. So we may take as dominating function
\begin{align*}
G_{a_z}(k) &= \frac{c^d}{2^d\pi^\frac{d}{2} }\left|\det\left(N_z^{-\frac{1}{2}}\cosh(zN)^{-1}\exp(zN)\right)\right| \exp\left(-\frac{c^2}{4} \sigma_z^\text{min} k^2\right)\\
\end{align*}
Where $\sigma_z^\text{min}$ denotes the lowest eigenvalue of $\left(\Re(N^{-1}_z) - M_p^2 \Re(N_z)^{-1}\right)$. For $G_{a_z}$ to be integrable, we require $\sigma_z^\text{min}>0$, or equivalently $\Re(N^{-1}_z) - M_p^2 \Re(N_z)^{-1}$ to be positive definite. As $\Re(N_z), \Im(N_z)$ are both in the functional calculus of the single self-adjoint operator $N$, they commute and so we find
\[\Re(N^{-1}_z) = \Re(N_z)(N_z^*N_z)^{-1}.\]
Using this, we find
\begin{equation}\label{eqnFuncCalcComp}
\Re(N^{-1}_z) - M_p^2 \Re(N_z)^{-1} = \left(\Re(N_z)^2(N_z^*N_z)^{-1} - M_p^2 I \right) \Re(N_z)^{-1}.
\end{equation}
But $\sigma(N_z)$ is a finite subset of $\Sigma_{\frac{\pi}{2}-\theta_p}$ so $\sigma(\Re(N_z)^{-1})$ is a finite subset of $(0, \infty)$, and so this is a product of commuting positive definite operators and is hence positive definite, noting $\Re(N_z)^2(N_z^*N_z)^{-1} = \cos^2(\text{arg}(N_z)) > M_p^2 I$. As integrability of $G_{a_z}$ implies boundedness of $a_z(Q,P)$, we find that $\exp(-zL)$ has a holomorphic extension as a $B(L^p(\mu))$-valued function from the domain $\mathbb{R}^+$ to the domain $E^N$.

Next we wish to investigate the semigroup for sets $E^N_{\epsilon, \delta} \subset E$, for all $\epsilon, \delta>0$ (see (\ref{eqnENdash})). For $z \in E^N_{\epsilon, \delta}$, we find that the dominating functions $G_{a_z}(k)$ are radially decaying for each $z$, and so the bound in Theorem \ref{thmRBoundedness0} becomes checking finiteness of:
\begin{align*}
&\sup_{z \in E^N_{\epsilon, \delta}} \int_{\mathbb{R}^d} G_{a_z}(k) dk \\
&= \sup_{z \in E^N_{\epsilon, \delta}} \int_{\mathbb{R}^d} \frac{c^d}{2^d\pi^\frac{d}{2} }\left|\det\left(N_z^{-\frac{1}{2}}\cosh(zN)^{-1}\exp(zN)\right)\right| \exp\left(-\frac{c^2}{4} \sigma_z^\text{min}k^2 \right) dk\\
&= \sup_{z \in E^N_{\epsilon, \delta}} \frac{c^d}{2^d\pi^\frac{d}{2} }\left|\det\left(N_z^{-\frac{1}{2}}\cosh(zN)^{-1}\exp(zN)\right)\right| \pi^\frac{d}{2}  \left(\frac{c^2}{4}  \sigma_z^\text{min} \right)^{-\frac{d}{2}}\\
&= \sup_{z \in E^N_{\epsilon, \delta}} \left|\det\left(\left(\sigma_z^\text{min}N_z\right)^{-\frac{1}{2}}\cosh(zN)^{-1}\exp(zN)\right)\right| \\
&= \sup_{z \in E^N_{\epsilon, \delta}} \left|\det\left(\sigma_z^\text{min}N_z\right)\right|^{-\frac{1}{2}} \left|\det\left(\cosh(zN)^{-1}\exp(zN)\right)\right|. 
\end{align*}
As the spectrum of $N$ is contained in the positive real line and $E^N_{\epsilon, \delta}$ is at least a distance of $\delta$ from the purely imaginary periodic poles of $\cosh(zN)^{-1}$, $\left|\cosh(zN)^{-1}\exp(zN)\right|$ is uniformly bounded for $z \in E^N_{\epsilon, \delta}$. Thus the second determinant above is uniformly bounded. As we have assumed $z \in E^N_{\epsilon, \delta}$, we find that $\left(\Re(N_z)^2(N_z^*N_z)^{-1} - M_p^2 I\right)  = \left(\cos^2(\text{arg}(N_z)) - M_p^2 I\right) > \epsilon I >0$, and so by Equation \ref{eqnFuncCalcComp}, $\sigma_z^\text{min}$ is bounded below by $\epsilon$ times the lowest eigenvalue of $(\Re(N_z))^{-1}$ which is the same as $\epsilon$ divided by the largest eigenvalue of $\Re(N_z)$. So expressing $\det\left(\sigma_z^\text{min}N_z\right)$ as the product of its eigenvalues, the supremum over $z \in E^N_{\epsilon, \delta}$ of $\left|\det\left(\sigma_z^\text{min}N_z\right)\right|^{-\frac{1}{2}}$ will be finite if and only if the ratio of largest eigenvalue of $\Re(N_z)$ and smallest eigenvalue of $N_z$ is uniformly bounded. We find by the spectral mapping theorem, that this ratio of eigenvalues is of the form
\[\frac{m}{n}\frac{\Re(\tanh(zn))}{|\tanh(zm)|},\]
where $n, m \in \sigma(N)$. As we have assumed $N$ is positive definite, $n, m>0$. Thus we find
\begin{align*}
\left| \frac{m}{n}\frac{\Re(\tanh(zn))}{|\tanh(zm)|} \right| &\leq  \frac{m}{n}\left|\frac{\tanh(zn)}{\tanh(zm)}\right| \\
&= \frac{m}{n}\left|\frac{(e^{2zn}-1)(e^{2zm}+1)}{(e^{2zn}+1)(e^{2zm}-1)}\right|.
\end{align*}
We denote this final bound $C_{n,m}(z)$. As $|\Re(z)| \to \infty$, $C_{n,m}(z)$ converges to $m/n$ uniformly in $\Im(z)$, so we may restrict to a set $\overline{E^n_{\epsilon, \delta}} \cap \{|\Re(z)| \leq C\}$. On such a set, $z$ is uniformly distant from the zeroes of $(e^{2zn}+1)(e^{2zm}-1)$, except $z=0$. However,
\begin{align*}
\lim_{z \to 0} \frac{m}{n}\frac{(e^{2zn}-1)(e^{2zm}+1)}{(e^{2zn}+1)(e^{2zm}-1)} &= \frac{m}{n} \lim_{z \to 0} \frac{(e^{2zn}-1)}{(e^{2zm}-1)} \\
&= \frac{m}{n} \lim_{z \to 0} \frac{2n e^{2zn}}{2m e^{2zm}} \\
&= 1.
\end{align*}
So $C_{n,m}(z)$ is bounded near $0$. Away from zero and with bounded real part, $C_{n,m}(z)$ is the product of two functions periodic in $\Im(z)$, $m\left|\frac{(e^{2zm}+1)}{(e^{2zm}-1)}\right|$ and $\frac{1}{n}\left|\frac{(e^{2zn}-1)}{(e^{2zn}+1)}\right|$, whose poles $z$ remains distant from. Using periodicity, boundedness is equivalent to boundedness on a compact set for each periodic function individually. However, both are continuous, and thus bounded. Thus $C_{n,m}(z)$ is uniformly bounded on $E^N_{\epsilon, \delta}$, and thus so is the relevant ratio of eigenvalues. As there are only finitely many choices for $m, n \in \sigma(N)$, we find that the relevant ratio of eigenvalues is uniformly bounded.
Hence 
\[\sup_{z \in E^N_{\epsilon, \delta}} \int_{\mathbb{R}^d} G_{a_z}(k) dk < \infty,\]
and so we apply Theorem \ref{thmRBoundedness0} to deduce that $\{a_z(Q,P); z \in E^N_{\epsilon, \delta}\}$ is R-bounded on $L^p(\mu)$.
\end{proof}

\begin{rem}
Both the domain $E^{\frac{1}{2}I}$, and the union of all domains of the form of $E^{\frac{1}{2}I}_{\epsilon, \delta}$ are exactly the classical Epperson region, which is known to be the largest domain on which the classical Ornstein-Uhlenbeck semigroup has a bounded analytic extension on $L^p(\mu)$ (see for example, \cite{Epperson}). The set $E^N$ is thus an analogue of the Epperson region, for certain variants of the classical Ornstein-Uhlenbeck operator. By examining how $E^{\frac{1}{2}I}_{\epsilon, \delta}$ fill out $E^{\frac{1}{2}I}$ as $\epsilon, \delta \to 0$, it can be seen that the Ornstein-Uhlenbeck semigroup is R-bounded if and only if it is uniformly bounded. This implies that the angle of sectoriality and R-sectoriality of the classical Ornstein-Uhlenbeck operator agree.
\end{rem}

\section{The Symbol Class $HS(M)$}\label{secMoreHSM}\label{sectionMoreHSM}

The first thing we wish to do is enrich the symbol class $HS_0(M)$ with an identity. In fact, without much more effort we can easily include symbols corresponding to anything in the Borel functional calculus of the position operators $Q$.

\begin{defn}
Let $B \subset L^\infty(\mathbb{R}^{2d})$ be the sub-Banach space $\{b \in L^\infty(\mathbb{R}^{2d}); b(x,\xi) = b(x,0), \text{ for a.e. } x,\xi \in \mathbb{R}^d\}$.
\end{defn}

\begin{lemma}\label{lemmaDecomposition}
$HS_0(M) \cap B = \{0\}$.
\end{lemma}

\begin{proof}
We will calculate $\mathcal{F}(I_b(x))$ for $b \in B$. Fixing some $\varphi \in \mathcal{S}(\mathbb{R}^d)$, we have:
\begin{align*}
\left\langle \mathcal{F}(I_b(x)), \varphi \right\rangle &= \left\langle I_b(x), \mathcal{F}^*\varphi \right\rangle \\
&= (2\pi)^{-\frac{d}{2}}\int_{\mathbb{R}^{2d}} b(x, \xi) \varphi(k)\exp(i\xi k)d\xi dk \\
&= (2\pi)^{-\frac{d}{2}} b(x,0) \int_{\mathbb{R}^{2d}} \varphi(k)\exp(i\xi k)d\xi dk \\
&= (2\pi)^{\frac{d}{2}} b(x,0) \varphi(0) \\
&= \left\langle (2\pi)^{\frac{d}{2}} b(x,0)\delta_0, \varphi \right\rangle.
\end{align*}
Where $\delta_0$ is the Dirac distribution. The second last equality follows from noting that $(2\pi)^{-d}\int_{\mathbb{R}^{2d}} \varphi(k)\exp(i\xi k)d\xi dk$ is the evaluation at $0$ of the inverse Fourier transform of the Fourier transform of $\varphi$. This clearly shows that $b$ does not satisfy the requirements of admission into $HS_0(M)$ unless $b=0$.
\end{proof}

\begin{defn}\label{defnHSM}
Fix $M:\mathbb{R}^d \to \mathbb{R}^d$. The space $HS(M)$ (standing for Holomorphic Strip) is a subspace of $L^\infty(\mathbb{R}^{2d})$, with
\[ HS(M) =  HS_0(M) \oplus B\]
with the norm of $a = a_0 + a_b$, $a_0 \in HS_0(M), a_b \in B$ defined by
\[||a||_{HS(M)} = ||a_0||_{HS_0(M)}+||a_b||_{L^\infty(\mathbb{R}^{2d})}.\]
(Note that this norm is well-defined as $HS_0(M) \cap B = \{0\}$, as in Lemma \ref{lemmaDecomposition}).
\end{defn}

Due to the generality of symbols in $B$, we can no longer use the formula of Theorem \ref{thmKernel} as a definition for the operator associated with a symbol in $B$. We thus provide an explicit extension of the generalised Weyl calculus to $B$ (as a contraction with respect to the $HS(M)$ norm and operator norm), motivated by our intuition as to how things should work. With this definition, the generalised Weyl calculus will extend to a bounded map $HS(M) \mapsto B(L^p(\mu))$.

\begin{defn}\label{defnMultiplierExt}
Define the extension of the generalised Weyl calculus to $a = a_b \in B$ via the action
\[(a(Q,P)f)(y) = a_b(y/c)f(y).\]
\end{defn}

This is natural because, formally, for $a = a_b \in B$ and $f \in C_\phi$ (so $\exp\left(-\frac{1}{2}\phi(\cdot)\right)f(\cdot) \in \mathcal{S}(\mathbb{R}^d)$), we have
\begin{align*}
(a(Q,P)f)(y) &= \frac{1}{(2\pi)^{d}c^d} \int_{\mathbb{R}^{2d}}a\left(\frac{x+y}{2c}, \xi\right) \exp\left(-i\xi \left(\frac{x-y}{c}\right)\right) \exp\left(\frac{1}{2}(\phi(y)+\phi(x))\right)f(x) d\xi d\mu(x) \\
&= \frac{1}{c^d}\exp\left(\frac{1}{2}\phi(y)\right)\int_{\mathbb{R}^{d}}a_b\left(\frac{x+y}{2c}\right)\delta_0\left(\frac{x-y}{c}\right)\exp\left(-\frac{1}{2}\phi(x)\right)f(x)dx \\
&= \frac{1}{c^d}\exp\left(\frac{1}{2}\phi(y)\right)\left\langle a_b\left(\frac{\cdot+y}{2c}\right)\delta_0\left(\frac{\cdot-y}{c}\right), \exp\left(-\frac{1}{2}\phi(\cdot)\right)f(\cdot)\right\rangle \\
&= a_b\left(\frac{2y}{2c}\right)\exp\left(\frac{1}{2}(\phi(y)-\phi(y))\right)f(y) \\
&= a_b\left(y/c\right)f(y).
\end{align*}
Where $\delta_0$ is the Dirac distribution. This extension is clearly contractive as a map $B \to B(L^p(\mu))$. Combining this with Theorem \ref{thmBoundedness0}, we have:

\begin{thm}\label{thmBoundednessComplete}
Fix $\phi \in C^2(\mathbb{R}^d)$ and $p \in [1, \infty]$, and suppose there exists a valid growth pair $(M,\epsilon)$ for $\phi$ and $p$.
Then the generalised Weyl calculus extends uniquely to a linear map $HS(M) \to B(L^p(\mu))$ and we have
\[||a(Q,P)||_{B(L^p(\mu))} \leq e^\epsilon||a||_{HS(M)}.\]
\end{thm}

Exactly as we have shown for $HS_0(M)$, we also get for free that symbols in $HS(M)$ have some holomorphic nature.

\begin{rem}[Holomorphic Nature of $HS(M)$]
Exactly as in Remark \ref{remHolomorphicExtension}, for any symbol $a \in HS(M)$ and $x \in \mathbb{R}^d$, $a(x, \xi+i\eta)$ has an extension as a function of $\xi + i\eta$ for $\eta = tM(x), t \in (-1,1)$, and the essential range of $a(x, \xi + i\eta)$ on domain $\{(x, \xi+i\eta)\in \mathbb{R}^d \times \mathbb{C}^d; \exists t \in (-1,1) \text{ s.t. } \eta = tM(x)\}$ is bounded by $||a||_{HS(M)}$. This follows trivially, as $a(x, \xi) = a_0(x, \xi)+a_b(x)$ and $a_0 \in HS_0(M)$ has the given extendability, while the constant function $\xi \mapsto a_b(x)$ has an entire and bounded (and constant) extension.
\end{rem}

A similar R-boundedness theorem also holds.

\begin{thm}\label{thmRBoundedness}
Fix $\phi \in C^2(\mathbb{R}^d)$ and $p \in [1, \infty]$, and suppose there exists a valid growth pair $(M,\epsilon)$ for $\phi$ and $p$.
Let $A \subset HS(M)$, and let $G \subset \mathcal{D}_{\infty,1}$, $A_B \subset B$ be the sets of $g_a, a_b$ corresponding to each $a \in A$ as in the decomposition $HS(M) = HS_0(M) \oplus B$ in Definition \ref{defnHSM}. Suppose that for each $g_a \in G$ we can choose dominating $G_a \in L^1(\mathbb{R}^d)$ such that the supremum over our selections of the quantity
\[\int_{\mathbb{R}^d} \esssup_{|y| \geq |x|} |G_a(y)| dx\]
is finite.
Also suppose that
\[\sup_{a_b \in B} ||a_b||_{L^\infty(\mathbb{R}^d)} < \infty.\] 
Then $A(Q,P) = \{a(Q,P), a \in A\}$ is R-bounded on $L^p(\mu)$.
\end{thm}

\begin{proof}
First note that $A \subset A_0 + A_B$, where $A_0$ and $A_B$ are the projections of $A$ onto $HS_0(M)$ and $B$ respectively, so by subadditivity of R-boundedness it suffices to check that $A_0(Q,P)$ and $A_B(Q,P)$ are R-bounded on $L^p(\mu)$. That $A_0(Q,P)$ is R-bounded follows from Theorem \ref{thmRBoundedness0}. Since everything in $A_B$ is in $B$, Definition \ref{defnMultiplierExt} implies that $A_B(Q,P)$ consists of multiplication operators $f(x) \mapsto a_b(x/c)f(x)$. Then since $\sup_{a_b \in B} ||a_b||_{L^\infty(\mathbb{R}^d)} < \infty$, $B(Q,P)$ is R-bounded on $L^q(\mu)$ for all $q \in (1, \infty)$, and hence on $L^p(\mu)$ (see, for example, Example 8.1.9 of \cite{AnalysisInBanachSpacesV2}). 
\end{proof}

We wish to study the space $HS(M)$ - equipped with a natural product (to be defined below) - to gain knowledge about operators on $L^p(\mu)$ related to $Q$ and $P$. To do so, we need to verify some facts about $HS(M)$.

\begin{thm}\label{thmHSMComplete}
For any measurable $M:\mathbb{R}^d \to \mathbb{R}^d$, $HS(M)$ is complete.
\end{thm}

\begin{proof}
Note that $HS(M)$ is a direct sum of the spaces $HS_0(M)$ and $B$, so providing both of these spaces are complete, we will be done. That $B$ is complete is obvious.
Let $\{a_n\} \subset HS_0(M)$ be a Cauchy sequence, with corresponding sequence $\{g_n\} \subset \mathcal{D}_{\infty,1}$

Then since $a_n$ is $HS_0(M)$-Cauchy, we find $\{g_n\}$ is $\mathcal{D}_{\infty,1}$-Cauchy, and hence has a limit $g \in \mathcal{D}_{\infty,1}$, say. Let $a \in L^\infty(\mathbb{R}^{2d})$ be for each $x, \xi \in \mathbb{R}^d$,
\[a(x, \xi) = (2\pi)^{-\frac{d}{2}}\int_{\mathbb{R}^d} (2\pi)^\frac{d}{2} g(x,k) \exp\left(-|M(x)k|\right)\exp(i k \xi) dk\]
Note this is well-defined for each $x, \xi$ as $\left|\exp\left(-|M(x)k|\right)\exp(i k \xi)\right| \leq 1$ and $|g(x,k)| < G(k)$ for some $G \in L^1(\mathbb{R}^d)$, and hence $|a(x,\xi)|$ is bounded by $||G||_{L^1(\mathbb{R}^d)}$. By the Fourier inversion formula on $\mathcal{S}(\mathbb{R}^d)$, $\mathcal{F}(I_a(x))\in \mathcal{S}'(\mathbb{R}^d)$ is given by integration against $(2\pi)^\frac{d}{2} g(x,k) \exp\left(-|M(x)k|\right)$, and so $a \in HS_0(M)$. It is clear that $||a_n-a||_{HS_0(M)} \to 0$ by construction. So $HS_0(M)$ is complete, and hence $HS(M)$ is complete.
\end{proof}

In the study of the standard Weyl calculus, there is a bilinear product $\mathcal{S}(\mathbb{R}^{2d}) \times \mathcal{S}(\mathbb{R}^{2d}) \to \mathcal{S}(\mathbb{R}^{2d})$, known as the Moyal product and denoted $\#$, which makes the Weyl calculus into an algebra homomorphism, I.e. such that for all $a_1, a_2 \in \mathcal{S}(\mathbb{R}^{2d})$, $a_1(X,D)a_2(X,D) = (a_1\#a_2)(X,D)$. We wish to define a similar  product on $HS(M)$, making the functional calculus an algebra homomorphism. Note that as our generalised Weyl pairs $(Q,P)$ are unitarily equivalent to the standard pair $(X,D)$ on $L^2$, any such product should agree with the Moyal product on $HS(M) \cap \mathcal{S}(\mathbb{R}^{2d})$, and hence we will also refer to such a product on $HS(M)$ as the Moyal product and denote it $\#$.

In the classical case, the Moyal product is either written in terms of an oscillatory integral involving $a_1(x, \xi)$ and $a_2(x, \xi)$, or as an asymptotic formula involving derivatives of $a_1$ and $a_2$. We will avoid both of these expressions by deducing what the product must be from the relation $a_1(Q,P)a_2(Q,P) = (a_1\#a_2)(Q,P)$. As both the definition of $HS(M)$ and the generalised Weyl calculus for $(Q,P)$ are written in terms of $\mathcal{F}_2a$ instead of $a$ explicitly, we will find a formula for the Moyal product in terms of $\mathcal{F}_2$ of the symbols.

\begin{thm}\label{thmMoyal}
For $a^1 = a^1_0+a^1_b, a^2 = a^2_0 + a^2_b \in HS(M) = HS_0(M) \oplus B$, the Moyal product $a^1\#a^2 \in L^\infty(\mathbb{R}^{2d})$ is given as the sum of $a^1_0\#a^2_b \in HS_0(M)$, $a^1_b\#a^2_0 \in HS_0(M)$, $a^1_b\#a^2_b \in B$, and $a^1_0\#a^2_0 \in L^\infty(\mathbb{R}^{2d})$ where
\begin{enumerate}
\item $\mathcal{F}_2(a^1_0\#a^2_0) = \left(2\pi\right)^{-\frac{d}{2}}\int_{\mathbb{R}^d}\mathcal{F}_2a^1_0\left( x+ \frac{v-k}{2}, v\right)\mathcal{F}_2a^2_0\left(x + \frac{v}{2}, k-v\right)dv$.
\item $\mathcal{F}_2(a^1_0\#a^2_b)(x, k) = \mathcal{F}_2a^1_0(x, k)a^2_b\left(x+\frac{k}{2}\right)$.
\item $\mathcal{F}_2(a^1_b\#a^2_0)(x, k) = a^1_b\left(x - \frac{k}{2}\right)\mathcal{F}_2a^2_0(x, k)$.
\item $(a^1_b\#a^2_b)(x) = a^1_b(x)a^2_b(x)$.
\end{enumerate}
Furthermore, suppose $M:\mathbb{R}^d \to \mathbb{R}^d$ is an affine function with real-symmetric linear part. Then $a^1_0\#a^2_0 \in HS_0(M)$, and the Moyal product is a Banach algebra product on $HS(M)$, such that for all $a^1, a^2 \in HS(M)$
\[||a^1 \# a^2||_{HS(M)} \leq ||a^1||_{HS(M)}||a^2||_{HS(M)}.\]
\end{thm}

\begin{proof}
We will only prove the first and second formula, as the third follows in almost the same way as the second, and the fourth is apparent from Definition \ref{defnMultiplierExt}. Recall Definition \ref{defnWeylHS0M}, which states for $a \in HS_0(M)$,
\[(a(Q,P)f)(y) = \frac{1}{(2\pi)^{d/2}c^d} \int_{\mathbb{R}^d} \mathcal{F}_2 a\left(\frac{x+y}{2c}, \frac{x-y}{c}\right) \exp\left(\frac{1}{2}(\phi(y)-\phi(x))\right)f(x)dx.\]
Thus we find
\begin{align*}
\left(a^1_0(Q,P)a^2_0(Q,P)f\right)(y) &= \frac{1}{(2\pi)^{d}c^{2d}} \int_{\mathbb{R}^{2d}} \mathcal{F}_2 a^1_0\left(\frac{z+y}{2c}, \frac{z-y}{c}\right) \exp\left(\frac{1}{2}(\phi(y)-\phi(z))\right)\\
&\mathcal{F}_2a^2_0\left(\frac{x+z}{2c}, \frac{x-z}{c}\right) \exp\left(\frac{1}{2}(\phi(z)-\phi(x))\right)f(x)dxdz \\
&= \frac{1}{(2\pi)^{d/2}c^{d}} \int_{\mathbb{R}^{d}} \left(\frac{1}{(2\pi)^{d/2}c^{d}}\int_{\mathbb{R}^d} \mathcal{F}_2 a^1_0\left(\frac{z+y}{2c}, \frac{z-y}{c}\right) \mathcal{F}_2a^2_0\left(\frac{x+z}{2c}, \frac{x-z}{c}\right)dz \right) \\
&\exp\left(\frac{1}{2}(\phi(y)-\phi(x))\right)f(x)dx
\end{align*}
and by making a change of variables $v = (z-y)/c$ we find
\begin{align*}
& \frac{1}{(2\pi)^{d/2}c^{d}}\int_{\mathbb{R}^d} \mathcal{F}_2 a^1_0\left(\frac{z+y}{2c}, \frac{z-y}{c}\right) \mathcal{F}_2a^2_0\left(\frac{x+z}{2c}, \frac{x-z}{c}\right)dz \\
&= \frac{1}{(2\pi)^{d/2}}\int_{\mathbb{R}^d} \mathcal{F}_2 a^1_0\left(\frac{x+y}{2c} - \frac{1}{2}\frac{x-y}{c} +\frac{v}{2}, v\right) \mathcal{F}_2a^2_0\left(\frac{x+y}{2c} + \frac{v}{2}, \frac{x-y}{c} -v\right)dv. \\
\end{align*}
By comparing with Definition \ref{defnWeylHS0M}, it can be seen that this confirms the formula for $\mathcal{F}_2(a^1_0\#a^2_0)(x, k)$.

For the second term, we find
\[\left(a^1_0(Q,P)a^2_b(Q,P)f\right)(y) = \frac{1}{(2\pi)^{d/2}c^d} \int_{\mathbb{R}^d} \mathcal{F}_2 a^1_0\left(\frac{x+y}{2c}, \frac{x-y}{c}\right) a^2_b\left(\frac{x}{c}\right)\exp\left(\frac{1}{2}(\phi(y)-\phi(x))\right)f(x)dx\]
and
\begin{align*}
\mathcal{F}_2 a^1_0\left(\frac{x+y}{2c}, \frac{x-y}{c}\right) a^2_b\left(\frac{x}{c}\right) &= \mathcal{F}_2 a^1_0\left(\frac{x+y}{2c}, \frac{x-y}{c}\right) a^2_b\left(\frac{x+y}{2c} + \frac{1}{2}\frac{x-y}{c}\right),
\end{align*}
which, after comparison with Definition \ref{defnWeylHS0M}, confirms the formula for $\mathcal{F}_2(a^1_0\#a^2_b)(x, k)$.

Note that by boundedness of $a^1_b, a^2_b$, the three products besides $a^1_0\#a^2_0$ lie in the spaces as given above. To see that the formula for $\mathcal{F}_2(a^1_0\#a^2_0)$ implies $a^1_0\#a^2_0 \in L^\infty(\mathbb{R}^{2d})$, we have
\begin{align*}
&\left|\mathcal{F}_2(a^1_0\#a^2_0)(x, k)\right| \leq \left(2\pi\right)^{-\frac{d}{2}}\int_{\mathbb{R}^d}\left|\mathcal{F}_2a^1_0\left( x+ \frac{v-k}{2}, v\right)\mathcal{F}_2a^2_0\left(x + \frac{v}{2}, k-v\right)\right|dv \\
&= \left(2\pi\right)^{\frac{d}{2}}\int_{\mathbb{R}^d}\exp\left(-\left|M\left(x+ \frac{v-k}{2}\right)v\right|-\left|M\left(x + \frac{v}{2}\right)(k-v)\right|\right)\left|g_{a^1_0}\left(x+ \frac{v-k}{2},v\right)g_{a^2_0}\left(x + \frac{v}{2},k-v\right)\right|dv\\
&\leq \left(2\pi\right)^{\frac{d}{2}}\int_{\mathbb{R}^d}\left|g_{a^1_0}\left(x+ \frac{v-k}{2},v\right)g_{a^2_0}\left(x + \frac{v}{2},k-v\right)\right|dv\\
&\leq \left(2\pi\right)^{\frac{d}{2}}\int_{\mathbb{R}^d}G_{a^1_0}\left(v\right)G_{a^2_0}\left(k-v\right)dv,\\
\end{align*}
which is in $L^1(\mathbb{R}^d)$ as the convolution of two $L^1$ functions. Hence $\mathcal{F}_2(a^1_0\#a^2_0) \in \mathcal{D}_{\infty,1}$, and so $a^1_0\#a^2_0 \in L^\infty(\mathbb{R}^{2d})$ as the partial inverse Fourier transform of a function dominated by an $L^1$ function (as in the proof of Theorem \ref{thmHSMComplete}).

We now suppose that $M$ is affine with real-symmetric linear part, say $M(x) = \tilde{M}(x)+\ell$ where $\tilde{M}$ is linear and real-symmetric, and $\ell \in \mathbb{R}^d$. Then
\begin{align*}
&\left|M(x)k\right|-\left|M\left(x+ \frac{v-k}{2}\right)v\right|-\left|M\left(x + \frac{v}{2}\right)(k-v)\right| \\
&= \left|\tilde{M}(x)k + \ell k\right|-\left|\tilde{M}\left(x\right)v + \tilde{M}\left(\frac{v-k}{2}\right)v+ \ell v\right|-\left|\tilde{M}\left(x\right)k + \ell k + \tilde{M}\left(\frac{v}{2}\right)k - \tilde{M}\left(x\right)v - \tilde{M}\left(\frac{v}{2}\right)v -\ell v\right| \\
&= \left|\tilde{M}(x)k + \ell k\right|-\left|\tilde{M}\left(x\right)v + \tilde{M}\left(\frac{v-k}{2}\right)v+ \ell v\right|-\left|\tilde{M}\left(x\right)k + \ell k - \tilde{M}\left(x\right)v - \tilde{M}\left(\frac{v-k}{2}\right)v -\ell v\right|,
\end{align*}
which is of the form $\left|A+B\right|-\left|A\right|-\left|B\right|$ for $A,B \in \mathbb{R}^d$, and is hence less than or equal to $0$ by the triangle inequality. Thus we find

\begin{align*}
\left|\mathcal{F}_2(a^1_0\#a^2_0)(x, k)\right| &\leq \left(2\pi\right)^{\frac{d}{2}}\int_{\mathbb{R}^d}\exp(-|M(x)k|)\exp\left(|M(x)k|-\left|M\left(x+ \frac{v-k}{2}\right)v\right|-\left|M\left(x + \frac{v}{2}\right)(k-v)\right|\right)\\
&\left|g_{a^1_0}\left(x+ \frac{v-k}{2},v\right)g_{a^2_0}\left(x + \frac{v}{2},k-v\right)\right|dv\\
&\leq \left(2\pi\right)^{\frac{d}{2}}\exp(-|M(x)k|)\int_{\mathbb{R}^d}\left|g_{a^1_0}\left(x+ \frac{v-k}{2},v\right)g_{a^2_0}\left(x + \frac{v}{2},k-v\right)\right|dv\\
&\leq \left(2\pi\right)^{\frac{d}{2}}\exp(-|M(x)k|)\int_{\mathbb{R}^d}G_{a^1_0}\left(v\right)G_{a^2_0}\left(k-v\right)dv\\
&= \left(2\pi\right)^{\frac{d}{2}}\exp(-|M(x)k|)G_{a^1_0}*G_{a^2_0}(k)
\end{align*}

Where $*$ denotes the convolution. Note that $G_{a^1_0}*G_{a^2_0} \in L^1(\mathbb{R}^d)$ as $G_{a^1_0}, G_{a^2_0} \in L^1(\mathbb{R}^d)$, so $a^1_0\#a^2_0 \in HS_0(M)$. Further,   $||G_{a^1_0}*G_{a^2_0}||_{L^1(\mathbb{R}^d)} \leq ||G_{a^1_0}||_{L^1(\mathbb{R}^d)}||G_{a^2_0}||_{L^1(\mathbb{R}^d)}$ so $||a^1_0\#a^2_0||_{HS_0(M)} \leq ||a^1_0||_{HS_0(M)}||a^2_0||_{HS_0(M)}$.

We check the other terms of $a^1 \# a^2$.
\begin{align*}
\mathcal{F}_2(a^1_0\#a^2_b)(x, k) &= \mathcal{F}_2a^1_0(x, k)a^2_b\left(x+\frac{k}{2}\right) \\
&= (2\pi)^\frac{d}{2}a^2_b\left(x+\frac{k}{2}\right)g_1(x,k) \exp\left(-M(x)|k|\right),
\end{align*}
with $a^2_b\left(x+\frac{k}{2}\right)g_1(x,k) \in \mathcal{D}_{\infty,1}$ with norm bounded by $||g_1||_{\mathcal{D}_{\infty,1}}||a^2_b||_{L^\infty(\mathbb{R}^{2d})}$.
\begin{align*}
\mathcal{F}_2(a^1_b\#a^2_0)(x, k) &= a^1_b(x - \frac{k}{2})\mathcal{F}_2a^2_0(x, k) \\
&= (2\pi)^\frac{d}{2}a^1_b\left(x-\frac{k}{2}\right)g_2(x,k) \exp\left(-M(x)|k|\right),
\end{align*}
with $a^1_b\left(x-\frac{k}{2}\right)g_2(x,k) \in \mathcal{D}_{\infty,1}$ with norm bounded by $||a^1_b||_{L^\infty(\mathbb{R}^{2d})}||g_2||_{\mathcal{D}_{\infty,1}}$.
\begin{align*}
(a^1_b\#a^2_0)(x) = a^1_b(x)a^2_b(x),
\end{align*}
with $a^1_b(x)a^2_b(x) \in B$ with $L^\infty$ norm bounded by $||a^1_b||_{L^\infty(\mathbb{R}^{2d})}||a^2_b||_{L^\infty(\mathbb{R}^{2d})}$.
Putting these together and using subadditivity of the $HS_0(M)$ norm gives
\begin{align*}
||a^1 \# a^2||_{HS(M)} &= ||a^1_0\#a^2_0+a^1_0\#a^2_b+a^1_b\#a^2_0||_{HS_0(M)} + ||a^1_0\#a^2_0||_{L^\infty(\mathbb{R}^{2d})}\\
&\leq ||a^1_0\#a^2_0||_{HS_0(M)}+||a^1_0\#a^2_b||_{HS_0(M)}+||a^1_b\#a^2_0||_{HS_0(M)} + ||a^1_0\#a^2_0||_{L^\infty(\mathbb{R}^{2d})}\\
&\leq ||a^1_0||_{HS_0(M)}||a^2_0||_{HS_0(M)}+||a^1_0||_{HS_0(M)}||a^2_b||_{L^\infty(\mathbb{R})}\\
&+||a^1_b||_{L^\infty(\mathbb{R})}||a^2_0||_{HS_0(M)} + ||a^1_0||_{L^\infty(\mathbb{R})}||a^2_0||_{L^\infty(\mathbb{R}^{2d})}\\
&=  \left(||a^1_0||_{HS_0(M)} + ||a^1_b||_{L^\infty(\mathbb{R})}\right)\left(||a^2_0||_{HS_0(M)} + ||a^2_0||_{L^\infty(\mathbb{R}^{2d})}\right)\\
&= ||a^1||_{HS(M)}||a^2||_{HS(M)}
\end{align*}
\end{proof}

We have the easily verified lemma and corollary:

\begin{lemma}
Assuming $M$ is part of a valid growth pair as in Theorem \ref{thmMoyal}, the Moyal product is associative, and has as identity the constant function $1(x, \xi) = 1$.
\end{lemma}

\begin{cor}\label{corBanachAlgebra}
Fix $\phi \in C^2(\mathbb{R}^d)$ and $p \in [1, \infty]$, and suppose there exists a valid growth pair $(M,\epsilon)$ for $\phi$ and $p$, with $M$ affine with real-symmetric linear part. Then $(HS(M), \#)$ is a unital Banach algebra, and the generalised Weyl calculus $HS(M) \to B(L^p(\mu))$, $a \mapsto a(Q,P)$ is a bounded Banach algebra homomorphism with norm at most $e^\epsilon$.
\end{cor}

This corollary makes our symbol class $HS(M)$ very distinct from the standard symbol classes of pseudodifferential calculus, and more like a single operator functional calculi. This suggests we really have the ``right" norm for symbols, or at least something very close. We can hypothesise that Corollary \ref{corBanachAlgebra} will allow us to get closer to bounded functional calculus for $L$ via softer Banach algebra techniques.

\section{Concluding Remarks}

\subsection{Semigroup Generation in $HS(M)$}\label{secSemigroupGeneration}

The application of the generalised Weyl calculus developed in Section \ref{secHInfty} was only possible because the symbol $a_t$ for the Ornstein-Uhlenbeck semigroup, such that $a_t(Q,P)=\exp(-tL)$, was known explicitly for $\phi$ of the form $\phi(x) = xN(x)+lx$. In this subsection we present some ideas about how one might be able to determine or prove existence of the symbol for the semigroup in cases other than such polynomials.

If we suppose that $a_t(Q,P) = \exp(-tL)$ for some symbol $a_t \in HS(M)$ and that we can extend the Moyal product to include products with $h$, of Theorem \ref{thmH} with $h(Q,P)=L$, then the semigroup properties of $\exp(-tL)$ correspond to the following properties of $a_t$
\begin{equation}\label{eqnHSMODE}\left\{\begin{array}{c}
\frac{da_t}{dt} = -h \# a_t, t > 0 \\
a_0 = 1.
\end{array}\right.
\end{equation}
We have replaced a (strong) ODE in the Banach algebra $B(L^p(\mu))$ with an ODE in the Banach algebra $HS(M)$. Since $HS(M)$ is also a space of functions, we have a lot of explicit tools to solve for the symbol $a_t$, such as taking ansatz. By taking a good ansatz, the above ODE can be solved explicitly when $\phi(x) = xN(x)+lx$, which is what was done by the author in determining the formula in Theorem \ref{TheoremSGHard}. This was also essentially the method used by H\"{o}rmander in \cite{UnpublishedHormander}, whose proof we refer to in the proof of our formula, although the method was discovered independently.

This method also lends itself to non-quadratic $\phi$, or perturbation, as an abstract semigroup generation problem in a Banach algebra. Noting Remark \ref{remAffine}, we know that any $\phi$ for which the theory presented in this paper is applicable must typically be a bounded $C^2$ perturbation of a potential of the form $xN(x)+lx$, and that the relevant $HS(M)$ class of symbols is the same for $\phi$ and $xN(x)+lx$. Thus it would be natural to consider the symbol for the semigroup of the relevant Ornstein-Uhlenbeck operator as a perturbation of the symbol for the semigroup for the Ornstein-Uhlenbeck associated with $xN(x)+lx$.

\subsection{Banach Algebra Techniques}

The idea of the $HS(M)$-valued ODE in Equation \ref{eqnHSMODE} could be taken even further. If we consider $h\#$ as an unbounded operator on the Banach space $HS(M)$ with an appropriate domain, we could ask when does $h\#$ have a bounded $H^\infty$ functional calculus, I.e. a bounded Banach algebra homomorphism $H^\infty(\Omega) \to B(HS(M))$, $f \mapsto f(h\#)$, for some domain $\Omega \subset \mathbb{C}$. Supposing such a bounded Banach algebra homomorphism exists, abstract theory of Banach algebras and the ``associative" nature of the unbounded operator $h\#$ (namely, that $h\#(a^1\#a^2)=(h\#a^1)\#a^2$) will imply that the image of the homomorphism will lie in the subspace of $B(HS(M))$ naturally identifiable with $HS(M)$ (as a Banach algebra $\mathcal{A}$ is always contained in $B(\mathcal{A})$). Making this identification, we can then compose with our generalised Weyl calculus, another bounded Banach algebra homomorphism, to obtain a bounded Banach algebra homomorphism $H^\infty(\Omega) \to B(L^p(\mu))$. As resolvents of $h\#$ would map to resolvents of $L$ under this homomorphism, it must be a bounded $H^\infty$ functional calculus for $L$ on $L^p(\mu)$. Diagrammatically, this process is as follows:

\begin{equation*}
\xymatrix@R=10pt{
    H^\infty(\Omega) \ar[r] & HS(M) \subset B(HS(M))\ar[r] & B(L^p) \\
    f \ar@{}[u]|{\rotatebox{90}{$\in$}} \ar@{|->}[r] 
            & f(h\#) \ar@{}[u]|{\rotatebox{90}{$\in$}} \ar@{|->}[r] 
            & f(h\#)(Q,P)=f(L) \ar@{}[u]|{\rotatebox{90}{$\in$}}
}
\end{equation*}

Thus we would find that the bounded $H^\infty$ functional calculus of $L$ factors through $HS(M)$, and so we would know that symbols for the functional calculus for $L$ exist in $HS(M)$, even if we can't write them explicitly. There are some interesting complications to this approach, such that the domain of $h\#$ cannot possibly be dense as its domain cannot include the subspace generated by the identity (or $B \subset HS(M)$, for that matter), while an assumption of dense domain is common in the literature of the $H^\infty$ functional calculus.

This leads us naturally to consider the spectrum of a symbol in $HS(M)$. We have seen in Remark \ref{remHolomorphicExtension} that symbols $a \in HS(M)$ have a ``pseudo-holomorphic" extension to the domain $D_M = \{(x, \xi+i\eta)\in \mathbb{R}^d \times \mathbb{C}^d; \exists t \in (-1,1) \text{ s.t. }\eta = t M(x)\}$, and that the image of this holomorphic extension is contained in the closed complex disc of radius $||a||_{HS(M)}$. It is true that the spectrum of an element of a Banach algebra is always contained in the set of complex numbers of modulus less than or equal to the norm of the element. Also, when we look at the multiplication operator subspace $B \subset HS(M)$, it is clear that the spectrum of an element is the essential range of the element, as an $L^\infty$ function of $x$. Similarly, if $M=0$ then $HS(M)$ will contain the Banach algebra $L^1(\mathbb{R}^d)$ with convolution as product (corresponding to those symbols which do not depend on $x$, in which case the Moyal product degenerates to convolution in $k$ of the partial Fourier transforms). It is true for $L^1(\mathbb{R}^d)$ with convolution as product, that the spectrum of an element is the range of its Fourier transform. Thus we could hypothesise that for $a \in HS(M)$, 
\[\sigma_{HS(M)}(a) := \{\lambda \in \mathbb{C}; (a-\lambda) \text{ is not invertible in } HS(M)\}\]
is related to the set
\[\text{EssRan}(a(x, \xi+i\eta)),\]
where the essential range is taken for the extension discussed in Remark \ref{remHolomorphicExtension}, over the domain $D_M = \{(x, \xi+i\eta)\in \mathbb{R}^d \times \mathbb{C}^d; \exists t \in (-1,1) \text{ s.t. }\eta = t M(x)\}$. 

While it might seem reasonable to think that the spectrum of a symbol is exactly the essential range of this extension of the symbol, that is most likely not true due to the following example. Consider one of the cases presented in Section \ref{secHInfty}, for $\phi(x) = \frac{x^2}{2}$, the unbounded operator $h\#$ with $h(x,\xi) = \frac{1}{2}(x^2+\xi^2-d)$ and $h(Q,P)=L$ the classical Ornstein-Uhlenbeck operator. In this case, a family of symbols for the semigroup generated by $L$ was found, which is uniformly bounded in $HS(M)$ for real time, and so $h\#$ generates (in some sense) a uniformly bounded semigroup in/on $HS(M)$ for real time. Thus $h\#$ should have spectrum with real part bounded below by $0$. However, the range of $h$ over the relevant domain for the analytic extension will always include the point $-\frac{d}{2}$, which should cause the semigroup to blow up for large time.

If there was some relationship between the spectrum of a symbol in $HS(M)$ and the range of its holomorphic extension, it would show that if $M \neq 0$, the only elements of $HS(M)$ with real spectrum are elements of $B$ which take real values. Assuming our argument about factorisation of the $H^\infty$ functional calculus through $HS(M)$ holds true, this could be seen as a more explicit reason as to why the classical Ornstein-Uhlenbeck operator has only holomorphic functional calculus on $L^p(\mu)$, for $p \neq 2$. Namely, the symbol $h(x, \xi) = \frac{1}{2}(x^2+\xi^2-d)$ for the classical Ornstein-Uhlenbeck operator will have non-real spectrum as an element of $HS(M)$, even though the classical Ornstein-Uhlenbeck operator has only real spectrum.

\bibliography{mybib}
\bibliographystyle{plain}

\end{document}